\documentclass[11pt, a4paper]{article}
\usepackage{amsmath}
\usepackage{amsfonts}
\usepackage{amssymb}
\usepackage{ulem}

\usepackage{amsthm}
\usepackage{url}
\usepackage{dcolumn}

\usepackage{cite}
\usepackage{fancyhdr}

\usepackage{graphicx, color}
\usepackage{hyperref}
\usepackage{xcolor}

\newtheorem*{rem*}{Remark}

\renewcommand{\thefigure}{\arabic{section}{-}\arabic{figure}}
\makeatletter
\@addtoreset{figure}{section}
\@addtoreset{table}{section}
\@addtoreset{table}{subsection}

\usepackage{here}

\begin{document}

\title{Converting non-periodic tilings with \mbox{Tile$(1, 1)$} into tilings with 
three types of pentagons, I}
\author{ Teruhisa SUGIMOTO$^{ 1), 2)}$ }
\date{}
\maketitle

{\footnotesize

\begin{center}
$^{1)}$ The Interdisciplinary Institute of Science, Technology and Art

$^{2)}$ Japan Tessellation Design Association

E-mail: ismsugi@gmail.com
\end{center}

}

{\small
\begin{abstract}
\noindent
Non-periodic tilings with \mbox{Tile$(1, 1)$} using the substitution method, 
as presented by Smith et al. in \cite{Smith_2024a} and \cite{Smith_2024b}, 
can be converted non-periodic tilings with three types of pentagons. 
When arbitrary replacements are excluded, the resulting non-periodic tilings with 
three types of pentagons exhibit two patterns. Note that, during the conversion process 
in this manuscript, the rhombus is not subdivided into smaller similar rhombuses.
\end{abstract}
}



\section{Introduction}
\label{section1}

\begin{figure}[t]
 \centering\includegraphics[width=15cm,clip]{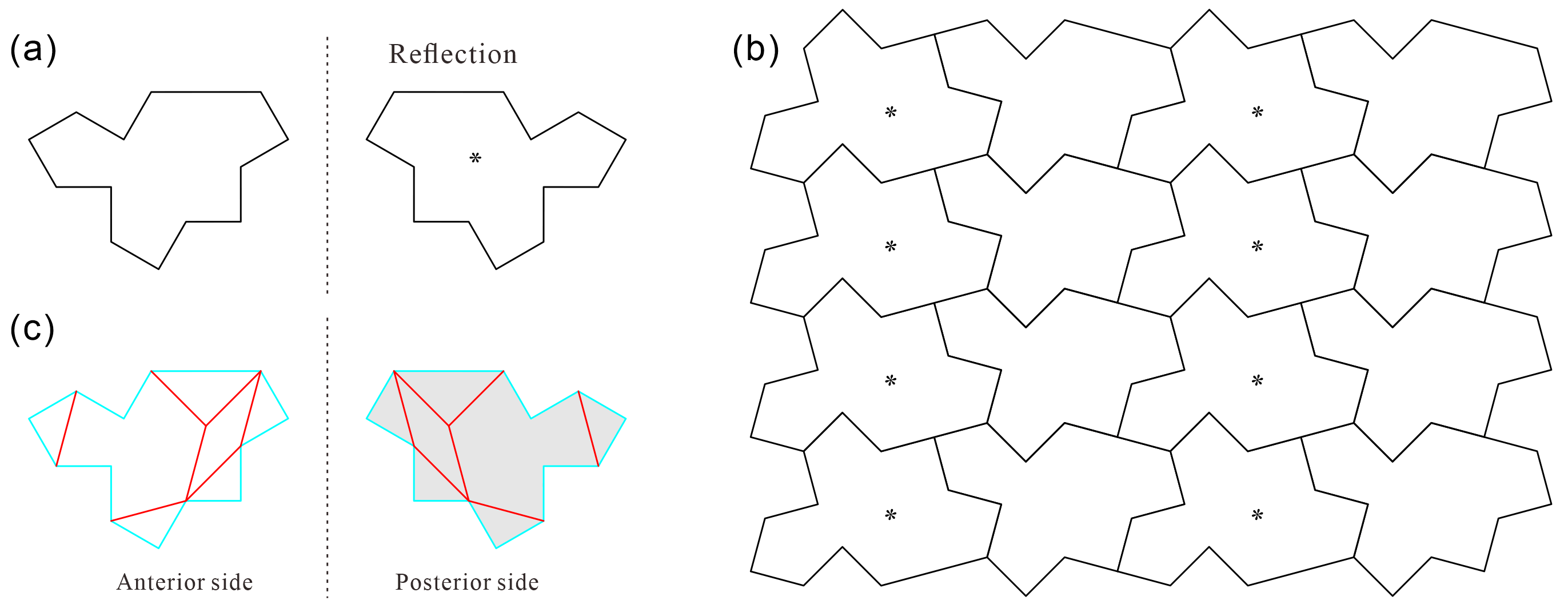} 
  \caption{{\small 
\mbox{Tile$(1, 1)$}, and periodic tiling formed by \mbox{Tile$(1, 1)$}.
} 
\label{Fig.1-1}
}
\end{figure}

Smith et al. \cite{Smith_2024a}. presented ``\mbox{Tile$(a, b)$}'' tiles and the substitution method, 
which can generate only non-periodic tiling\footnote{ 
A tiling (or \textit{tessellation}) of the plane is a collection of sets, called tiles, that cover the 
plane without gaps or overlaps, except for the boundaries of the tiles. The term ``tile'' refers 
to a topological disk, whose boundary is a simple closed curve \cite{G_and_S_1987}. A tiling 
exhibits \textit{periodicity} if its translation by a non-zero vector coincides with itself; a 
tiling is considered periodic if it coincides with its translation by two linearly independent vectors. 
However, in this study, a tiling with periodicity is referred to as \textit{periodic}, and a tiling 
without periodicity is referred to as \textit{non-periodic}.
}  All tiles in the tiling are of the same size and shape, but the tiling generation process requires 
the use of reflected tiles. Therefore, the tiles are referred to as ``aperiodic monotiles\footnote{ 
If all tiles in a tiling are of the same size and shape, the tiling is described as \textit{monohedral}, 
which allows the use of reflected tiles (the posterior side tiles) in monohedral tilings \cite{G_and_S_1987}. 
In other words, in monohedral tilings, the anterior and posterior sides of the tiles are treated as 
the same type (i.e., the concept assumes that there is only one type of tile). A tile that admits a 
monohedral tiling is referred to as a \textit{monotile}.
}.'' \mbox{Tile$(1, 1)$}, which corresponds to $a=b = 1$ in \mbox{Tile$(a, b)$}, can generate a periodic 
tiling if reflected tiles are allowed to be used during tiling generation (see Figure~\ref{Fig.1-1}(a) and 
(b)) \cite{Smith_2024a}. In a subsequent study, Smith et al. \cite{Smith_2024b} showed that 
\mbox{Tile$(1, 1)$} can generate only non-periodic tiling if its reflection is prohibited during tiling 
generation. This \mbox{Tile$(1, 1)$} is referred to as ``chiral aperiodic monotile,\footnote{ 
To be precise, \mbox{Tile$(1, 1)$} is considered a ``weakly chiral aperiodic monotile \cite{Smith_2024b}.''
}'' and the corresponding substitution method was described in \cite{Smith_2024b}.

\begin{figure}[htbp]
 \centering\includegraphics[width=15cm,clip]{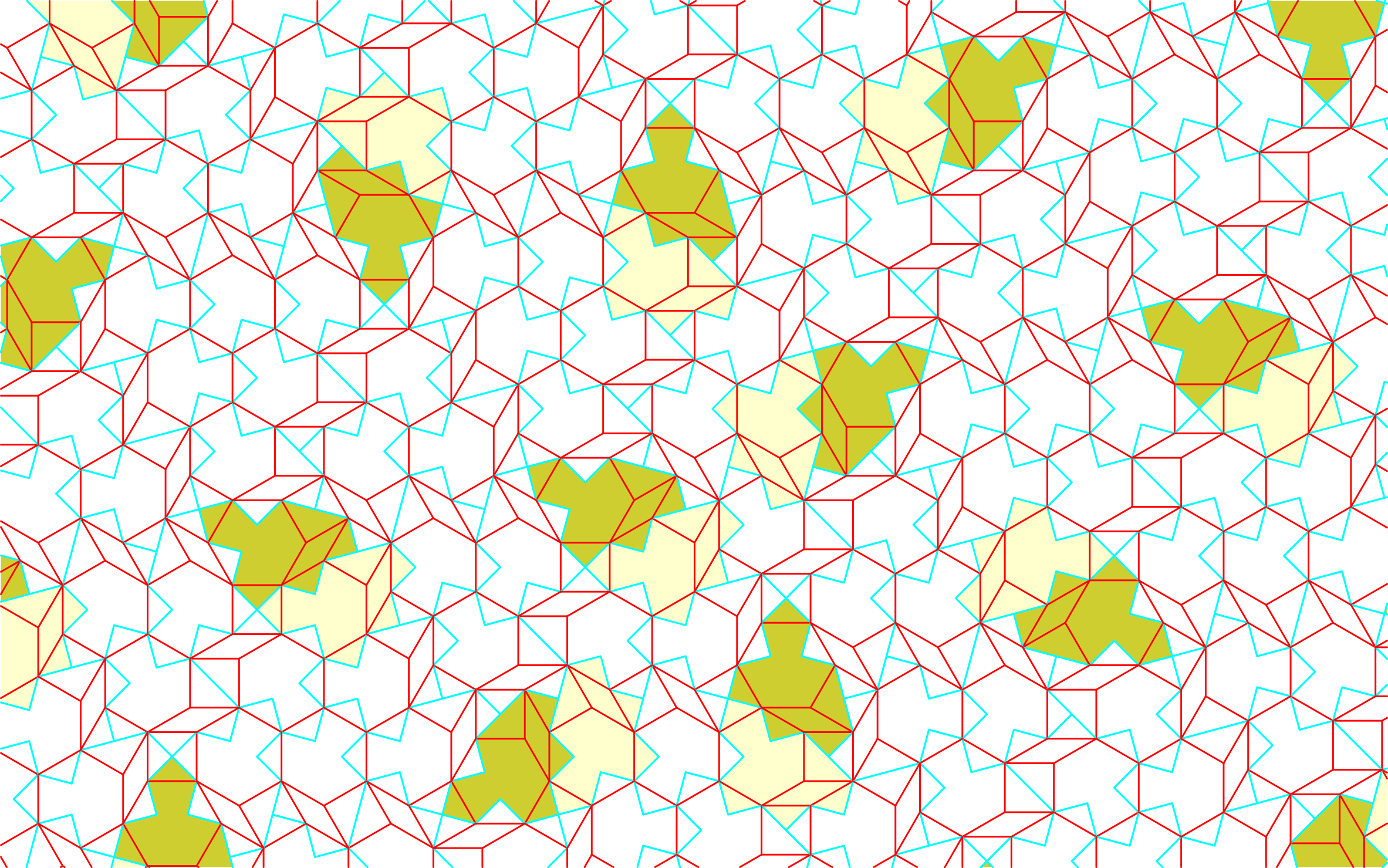} 
  \caption{{\small 
Converting non-periodic tiling $T_{s}$ with \mbox{Tile$(1, 1)$} into tiling consisting of 
squares, regular hexagons, and rhombuses with an acute angle of $30^ \circ$.
} 
\label{Fig.1-2}
}
\end{figure}

\begin{figure}[htbp]
 \centering\includegraphics[width=15cm,clip]{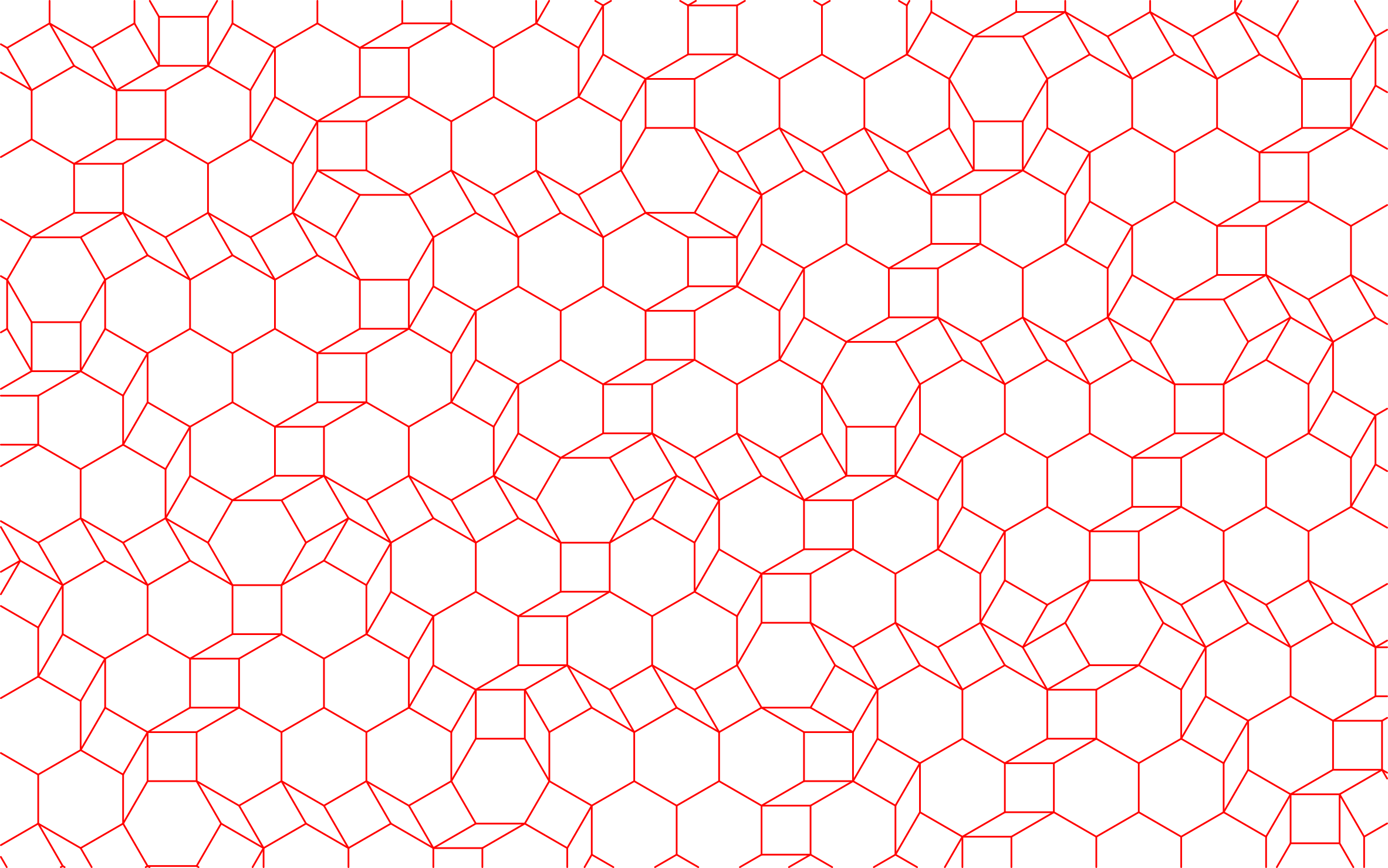} 
  \caption{{\small 
Non-periodic tiling consisting of squares, regular hexagons, and rhombuses with 
an acute angle of $30^ \circ$ generated by conversion based on 
$T_{s}$ with \mbox{Tile$(1, 1)$}.
} 
\label{Fig.1-3}
}
\end{figure}

\begin{figure}[htbp]
 \centering\includegraphics[width=15cm,clip]{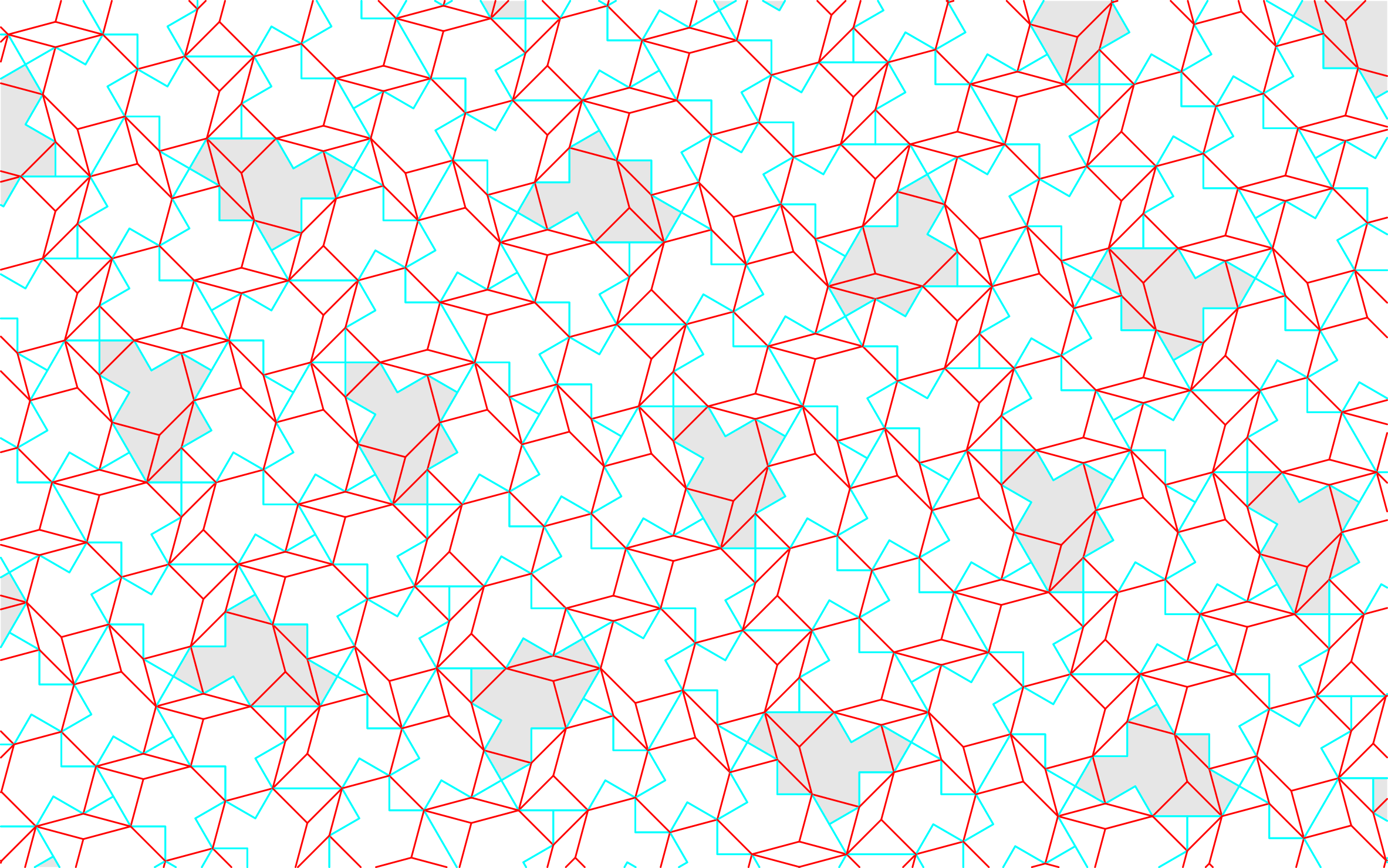} 
  \caption{{\small 
Converting non-periodic tiling $T_{h}$ with \mbox{Tile$(1, 1)$} into tiling consisting of 
squares, regular hexagons, and rhombuses with an acute angle of $30^ \circ$.
} 
\label{Fig.1-4}
}
\end{figure}

\begin{figure}[htbp]
 \centering\includegraphics[width=15cm,clip]{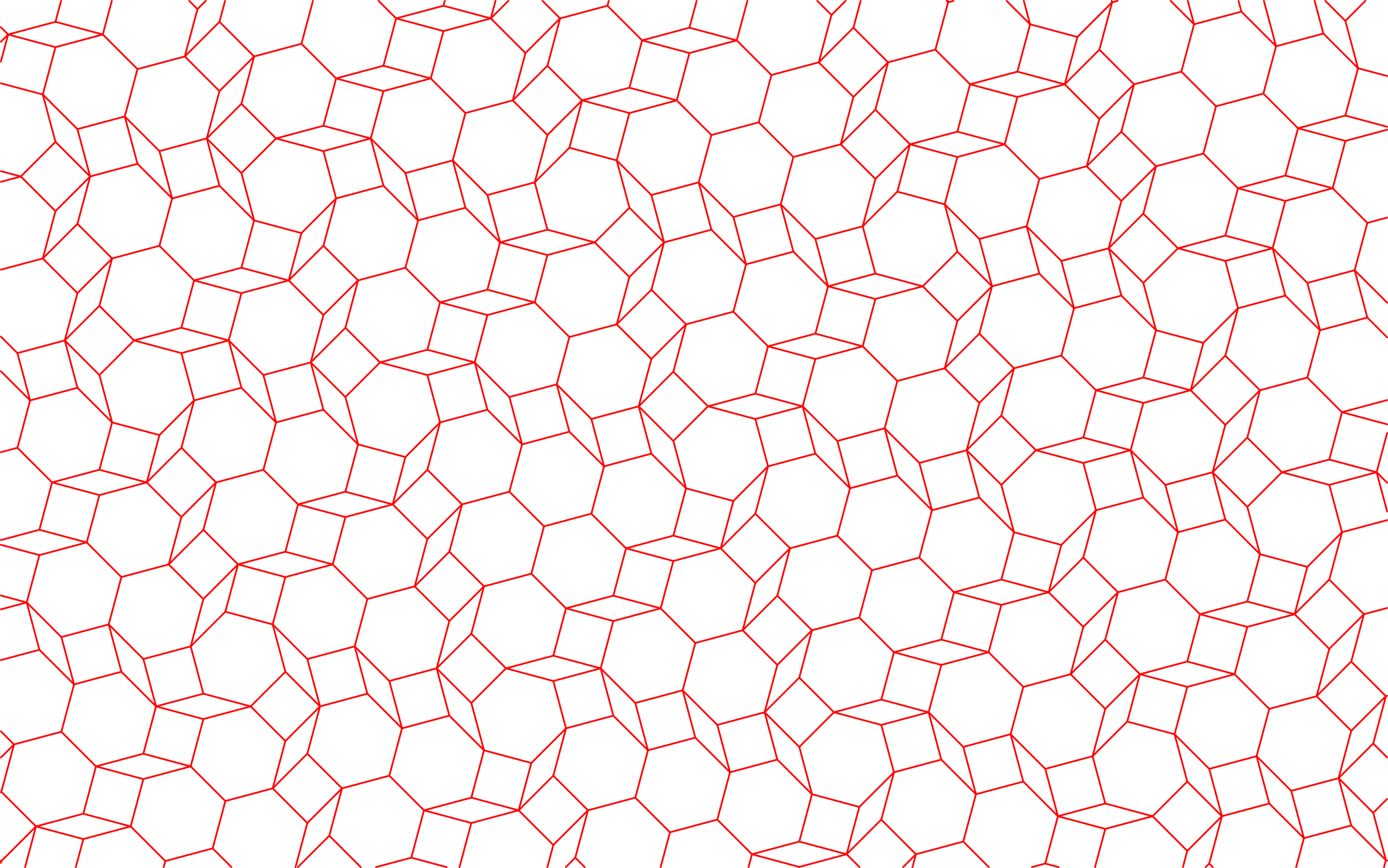} 
  \caption{{\small 
Non-periodic tiling consisting of squares, regular hexagons, and rhombuses with 
an acute angle of $30^ \circ$ generated by conversion based on
$T_{h}$ with \mbox{Tile$(1, 1)$}.
} 
\label{Fig.1-5}
}
\end{figure}

In this study, let $T_{h}$ be a non-periodic tiling that can be generated using the clusters 
$H_{7}$ and $H_{8}$ in conjunction with the substitution method, as described in \cite{Smith_2024a}. 
Furthermore, let $T_{s}$ be a non-periodic tiling that can be generated using the clusters in 
conjunction with the substitution method, as described in the Appendix of \cite{Smith_2024b}. 
If \mbox{Tile$(1, 1)$} is assigned a pattern incorporating decomposition lines, as shown in 
Figure~\ref{Fig.1-1}(c), then $T_{s}$ with \mbox{Tile$(1, 1)$} can be converted into a tiling 
consisting of squares, regular hexagons, and rhombuses with an acute angle of $30^ \circ$, 
as shown in Figures~\ref{Fig.1-2} and \ref{Fig.1-3}\footnote{ 
This conversion was reported in a mailing list of ``Tiling Lovers.''
}. Similarly, using \mbox{Tile$(1, 1)$} with the pattern incorporating decomposition lines, 
as shown in Figure~\ref{Fig.1-1}(c), $T_{h}$ with \mbox{Tile$(1, 1)$} can be converted into a 
tiling consisting of squares, regular hexagons, and rhombuses with an acute angle of $30^ \circ$ 
as in Figures~\ref{Fig.1-4} and \ref{Fig.1-5}. Because a regular hexagon can be divided into 
three rhombuses with an acute angle of $60^ \circ$, both $T_{h}$ and $T_{s}$ with 
\mbox{Tile$(1, 1)$} can be converted into non-periodic tilings consisting of rhombuses with 
acute angles of $90^ \circ$, $60^ \circ$, and $30^ \circ$. (Note: A rhombus with an acute 
angle of $90^ \circ$ is a square). However, the division of a regular hexagon (the orientation 
of rhombuses obtained by division) is not uniquely determined.

As demonstrated in \cite{Sugimoto_2020} and \cite{Sugimoto_2022}, a rhombic edge-to-edge 
tiling\footnote{ 
An edge-to-edge tiling with polygons is defined as a tiling in which the vertices (corners) and edges 
(sides) of the polygons coincide with the vertices (points where three or more tiles meet) and edges 
of the tiling \cite{G_and_S_1987}.
} (i.e., an edge-to-edge tiling with rhombuses) can be converted into a pentagonal tiling 
(i.e., a tiling with pentagons) using a specific method. We examined this conversion method and 
non-periodic tilings $T_{h}$ and $T_{s}$ consisting of rhombuses with acute angles of $90^ \circ$, 
$60^ \circ$, and $30^ \circ$ based on Figures~\ref{Fig.1-3} and \ref{Fig.1-5}. As a result, we 
confirmed that non-periodic tilings $T_{h}$ and $T_{s}$ with \mbox{Tile$(1, 1)$} can be converted 
into non-periodic tilings with three types of pentagons. Furthermore, when arbitrary replacements 
are excluded, we found that each non-periodic tiling with three types of pentagons, which is 
generated by conversion based on $T_{h}$ and $T_{s}$ with \mbox{Tile$(1, 1)$}, exhibits two 
patterns (design patterns created by the arrangement of polygonal tiles). This manuscript 
presents the results of converting rhombic edge-to-edge tilings into pentagonal tilings 
without performing a similarity division on the rhombus (i.e., subdivision of the rhombus into 
smaller similar rhombuses).

\renewcommand{\figurename}{{\small Figure}}

\renewcommand{\figurename}{{\small Figure}}

\renewcommand{\figurename}{{\small Figure}}

\renewcommand{\figurename}{{\small Figure}}

\renewcommand{\figurename}{{\small Figure}}



\section{Preparation}
\label{section2}

Before presenting the results of this study, we provide a brief review of 
the conversion of an edge-to-edge tiling with rhombuses into a tiling with 
pentagons, as presented in \cite{Sugimoto_2020} and \cite{Sugimoto_2022}.

\begin{figure}[htbp]
 \centering\includegraphics[width=15cm,clip]{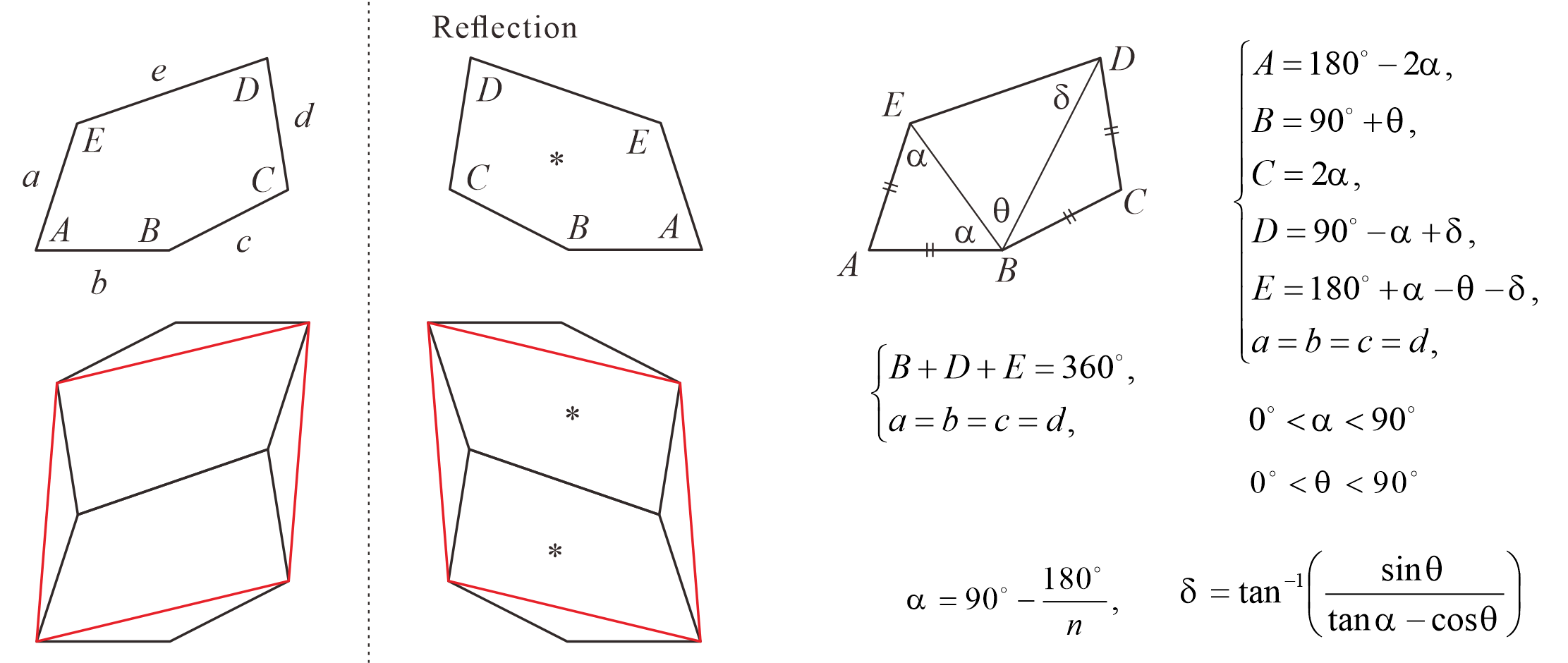} 
  \caption{{\small 
Properties of the pentagon used in the conversion, and the relationships between 
the rhombus and the pair of pentagons (pentagonal pair).
} 
\label{Fig.2-1}
}
\end{figure}

\begin{itemize}
	\item[(i)] An edge-to-edge tiling with rhombuses can be converted into a tiling with 
	pentagons by converting one rhombus into two pentagons (a pair of pentagons)\footnote{ 
	All quadrangles (rhombuses) are monotiles, and the pentagons used in this conversion 
	are also monotiles. If the pentagons used in this conversion are convex, then they 
	are classified as convex pentagonal monotiles belonging to the Type 2 family 
	\cite{G_and_S_1987, Sugimoto_2020, Sugimoto_2022}. The outline of this pair of 
	pentagons is regarded as a parallel equilateral concave octagon (a concave octagon 
	in which all four pairs of opposite edges are parallel and all edges are of equal length). 
	In other words, it can be regarded as a rhombus transformed into a parallel equilateral 
	concave octagon which can generate a tiling by translation.
	}, as shown in Figure~\ref{Fig.2-1}.
	\item[(ii)] Because a rhombus can be divided into its own similar figures, 
	a tiling with rhombuses can be converted into a tiling with pentagons after 
	the original rhombus is divided such that it contains $4 \cdot u^2$ for 
	$u = 1,\, 2,\, 3,\ldots $ of its own similar rhombuses.
	\item[(iii)] The conversion of a tiling consisting of rhombuses with acute angles of 
	$90^ \circ$, $60^ \circ$, and $30^ \circ$ into a pentagonal tiling can be performed 
	using the contents of Section 5.3 in \cite{Sugimoto_2022}. Depending on the value of 
	the parameter $\theta$, seven patterns (Sections 5.3.1 to 5.3.7 in \cite{Sugimoto_2022}) 
	with different combinations of three types of pentagons are realizable.
\end{itemize}

\begin{figure}[t]
 \centering\includegraphics[width=15cm,clip]{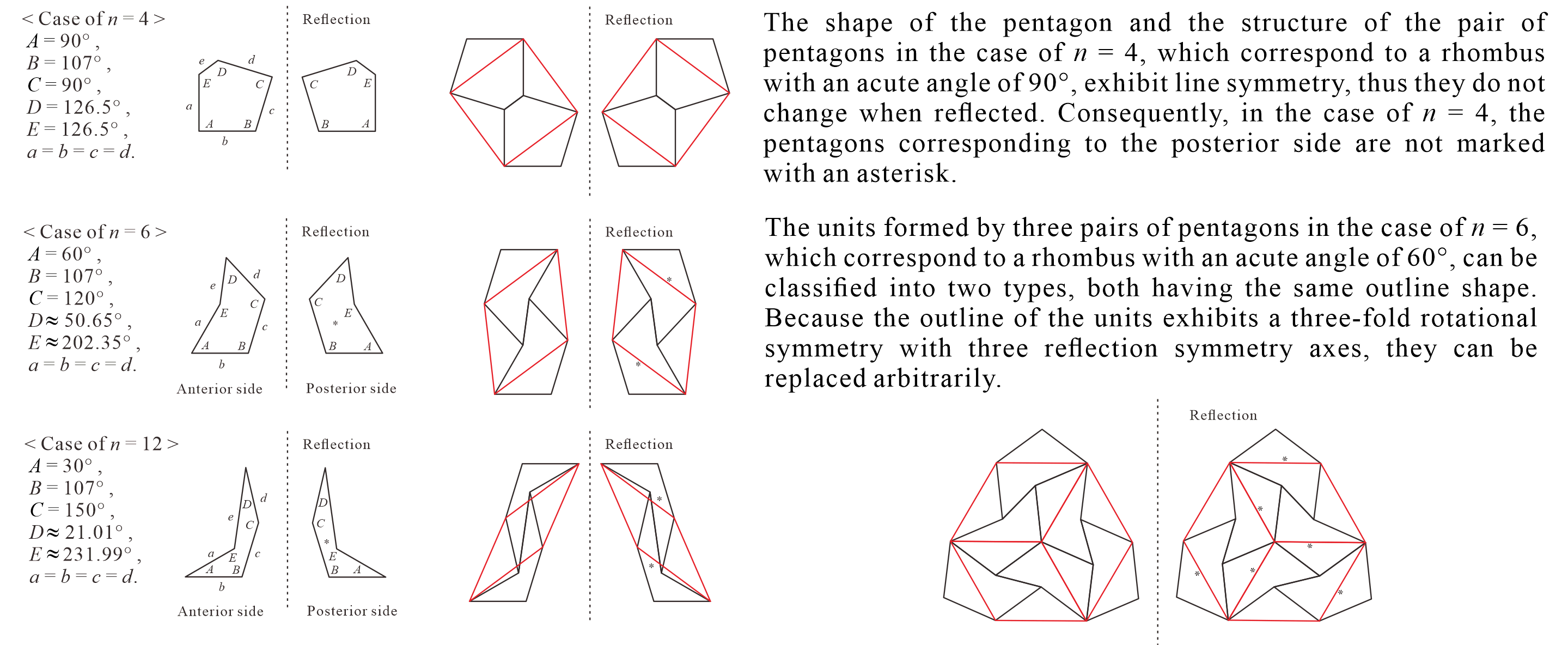} 
  \caption{{\small 
Pentagons used for the conversion in this study.
} 
\label{Fig.2-2}
}
\end{figure}

In this study, the three types of pentagons shown in Figure~\ref{Fig.2-2}, derived from 
Section 5.3.1 ($\theta = 17^ \circ $ in Table 5.3.1, i.e., $B = 107^ \circ$) of \cite{Sugimoto_2022}, 
were used for the conversions. In other words, the following conversions were applied:

\begin{itemize}
	\item Convert the rhombus with an acute angle of $90^ \circ$ into a pair 
	of convex pentagons.
	\item Convert the rhombuses with acute angles of $30^ \circ$ and 
	$60^ \circ$ into pairs of concave pentagons.
\end{itemize}

Therefore, the non-periodic tilings $T_{h}$ and $T_{s}$ are converted into tilings with one 
type of convex pentagon and two types of concave pentagons. Owing to the properties 
of the conversion method, conversions of other patterns that vary depending on the 
parameter $\theta$, such as conversions using three types of convex pentagons, 
can also be obtained in the same way as the conversions performed in this study 
\cite{Sugimoto_2022, Sugimoto_site_Ts, Sugimoto_site_Th}.

\renewcommand{\figurename}{{\small Figure}}

In this manuscript, the conversion presented in (ii), where a tiling with rhombuses is 
converted into a tiling with pentagons after the original rhombus is divided such that 
it contains its own similar rhombuses (i.e., the original rhombus is subdivided into 
smaller similar rhombuses), was not performed. Therefore, the results and discussions 
presented in this manuscript exclude the case in which the conversion is applied after dividing 
the original rhombus (i.e., a rhombus whose edge length is equal to the decomposition line 
assigned to \mbox{Tile$(1, 1)$} in Figure~\ref{Fig.1-1}(c)) into its own similar rhombuses.

The orientation of rhombuses with an acute angle of $60^ \circ$, obtained by dividing a 
regular hexagon into three sections, can be classified into two types and can be replaced 
arbitrarily. In the tiling of Figures~\ref{Fig.1-3} and \ref{Fig.1-5}, a rhombus with an acute 
angle of $60^ \circ$ appears if and only if a regular hexagon is divided into three; however, 
no rule in this conversion has been identified that determines the orientation of the 
rhombuses obtained by dividing a regular hexagon into three sections. In the pentagonal 
tilings obtained by conversion in this study, the arbitrary replacement corresponds to 
swapping a unit formed by three pairs of pentagons on the anterior side and a unit formed 
by three pairs of pentagons on the posterior side (i.e., reflection) (see Figure~\ref{Fig.2-2}). 
In this study, during the conversion in the case where the original rhombus is not subdivided 
into smaller similar rhombuses, the arbitrary replacement is not performed, and all pentagons 
corresponding to the regular hexagon are on the anterior side (i.e., pentagons marked without 
asterisks in Figure~\ref{Fig.2-2}). However, in \ref{appA}, we provide a brief introduction to 
the case in which this arbitrary replacement is applied.

\renewcommand{\figurename}{{\small Figure}}



\section{Conversion of non-periodic tiling $T_{h}$ in the case where the 
original rhombus is not subdivided into smaller similar rhombuses}
\label{section3}

\begin{figure}[htbp]
 \centering\includegraphics[width=14.5cm,clip]{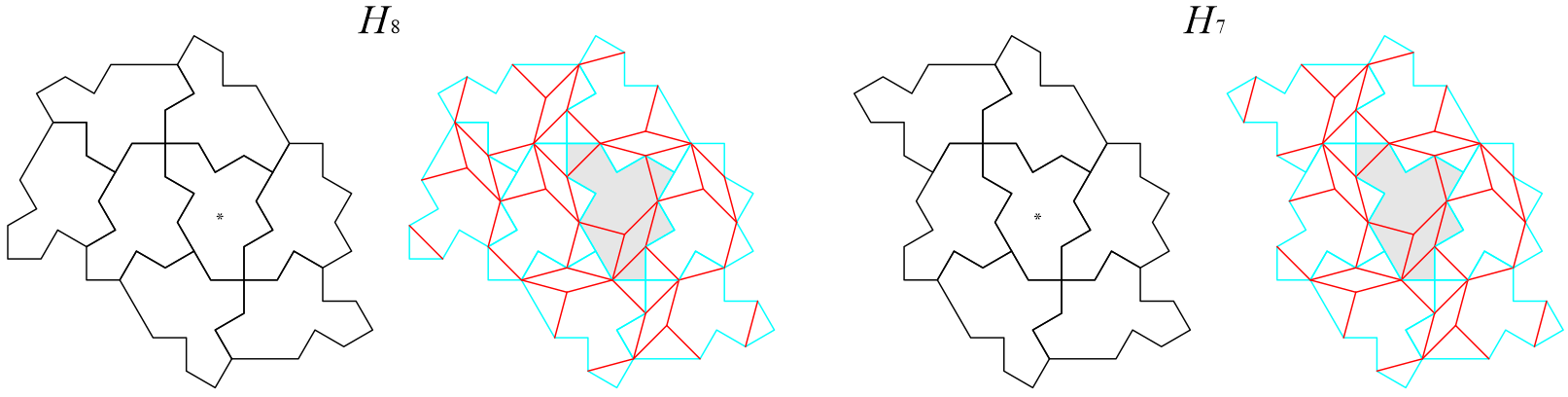} 
  \caption{{\small 
Clusters $H_{8}$ and $H_{7}$ formed by \mbox{Tile$(1, 1)$}.
} 
\label{Fig.3-1}
}
\end{figure}

When $T_{h}$ with \mbox{Tile$(1, 1)$} is converted into a tiling with three types of pentagons 
without subdivision of the original rhombus into smaller similar rhombuses (i.e., in the case where 
the original rhombus is not subdivided into smaller similar rhombuses), two patterns are observed: 
(i) all pentagons in the pairs corresponding to rhombuses with an acute angle of $30^ \circ$ 
are marked with asterisks (the posterior side in Figure~\ref{Fig.2-2}), or 
(ii) they are marked without asterisks (the anterior side in Figure~\ref{Fig.2-2}). 
Additionally, the patterns on \mbox{Tile$(1, 1)$} for forming each pattern of the pentagonal 
tilings (i.e., the patterns created by the decomposition lines corresponding to each pentagonal 
tiling) can be classified into multiple types. In other words, a single pattern such as that 
presented in Figure~\ref{Fig.1-1}(c) is not observed. In this study, we do not classify 
\mbox{Tile$(1, 1)$} based on different patterns, nor do we use different patterns of 
\mbox{Tile$(1, 1)$} to generate tilings. In this study, the conversion results were obtained by 
generating tilings using clusters $H_{8}$ and $H_{7}$, in conjunction with the substitution method 
described in \cite{Smith_2024a}. 
In other words, we generated non-periodic tilings with three types of pentagons using the 
clusters composed of three types of pentagons based on $H_{8}$ and $H_{7}$ formed by 
\mbox{Tile$(1, 1)$} (see Figure~\ref{Fig.3-1}) , in conjunction with the substitution method 
for $H_{8}$ and $H_{7}$ described in \cite{Smith_2024a}.

\renewcommand{\figurename}{{\small Figure}}


\subsection{First pattern of $T_{h}$ in the case where the original rhombus is not 
subdivided into smaller similar rhombuses}
\label{subsection3.1}

\begin{figure}[htbp]
 \centering\includegraphics[width=15cm,clip]{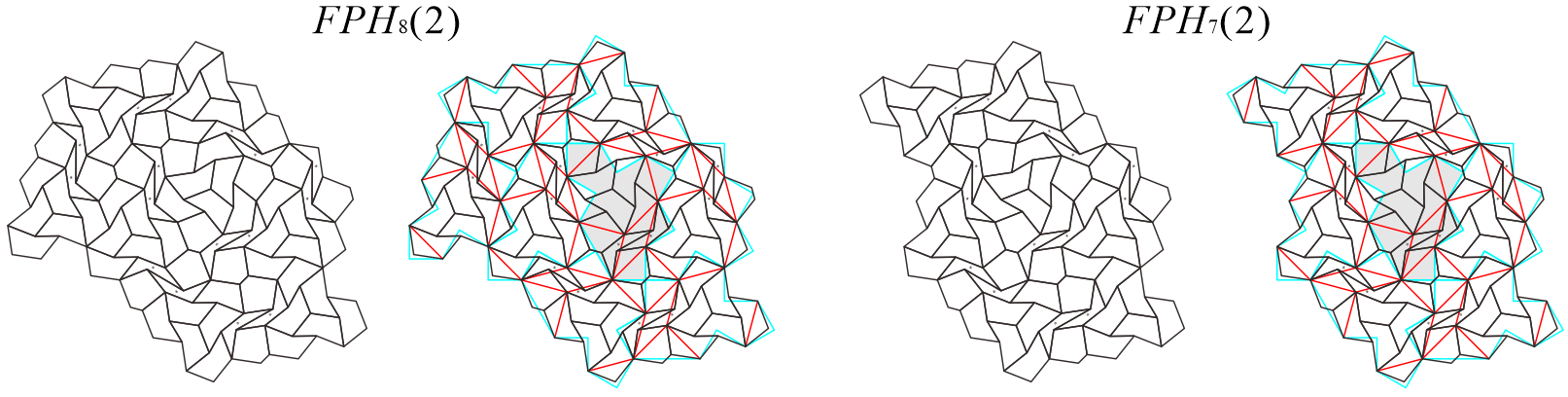} 
  \caption{{\small 
$FPH_{8}(2)$ and $FPH_{7}(2)$.
} 
\label{Fig.3-2}
}
\end{figure}

In the first pattern of $T_{h}$ (the first pattern when $T_{h}$ with \mbox{Tile$(1, 1)$} is converted 
into a tiling with three types of pentagons) in the case where the original rhombus is not subdivided 
into smaller similar rhombuses, all pentagons in the pairs corresponding to rhombuses with an 
acute angle of $30^ \circ$ are marked with asterisks (the posterior side in Figure~\ref{Fig.2-2}). 
The clusters formed by the three types of pentagons based on clusters $H_{8}$ and $H_{7}$ 
(see Figure~\ref{Fig.3-1}) for this pattern are shown in Figure~\ref{Fig.3-2}. Let these 
clusters be $FPH_{8}(2)$ and $FPH_{7}(2)$.

\begin{figure}[htbp]
 \centering\includegraphics[width=15cm,clip]{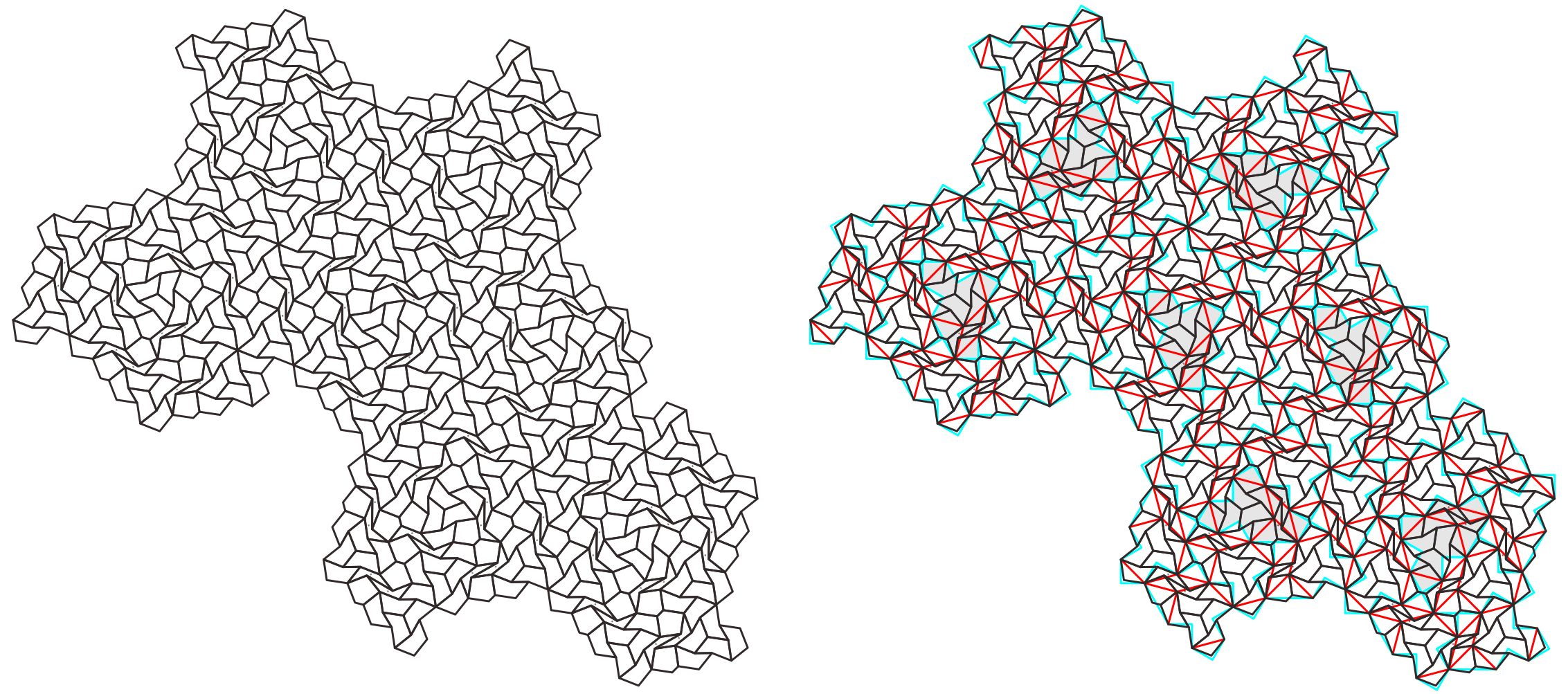} 
  \caption{{\small 
$FPH_{8}(3)$.
} 
\label{Fig.3-3}
}
\end{figure}

Figure~\ref{Fig.3-3} shows cluster $FPH_{8}(3)$ of the next substitution step formed by 
$FPH_{8}(2)$ and $FPH_{7}(2)$. Figure 3-4 shows the first pattern of $T_{h}$ with three types 
of pentagons in the case where the original rhombus is not subdivided into smaller similar rhombuses, 
which was generated using $FPH_{8}(2)$ and $FPH_{7}(2)$ in conjunction with the substitution method 
described in \cite{Smith_2024a}.

\renewcommand{\figurename}{{\small Figure}}

\renewcommand{\figurename}{{\small Figure}}

\vspace{1mm}

\renewcommand{\figurename}{{\small Figure}}
\begin{figure}[htbp]
 \centering\includegraphics[width=15cm,clip]{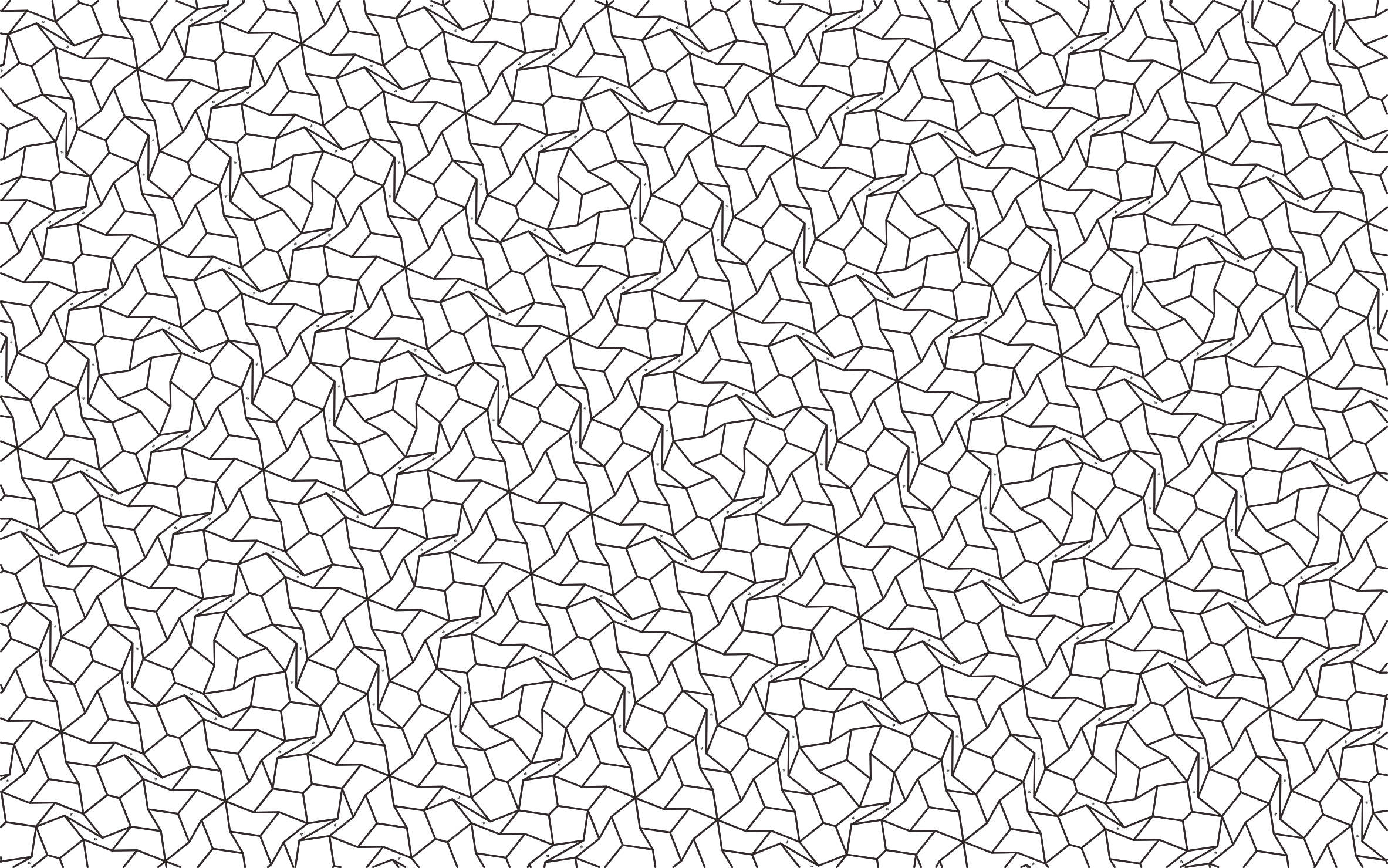} 
  \caption{{\small 
First pattern of $T_{h}$ with three types of pentagons in the case where the original 
rhombus is not subdivided into smaller similar rhombuses.
} 
\label{Fig.3-4}
}
\end{figure}


\subsection{Second pattern of $T_{h}$ in the case where the original rhombus is not 
subdivided into smaller similar rhombuses}
\label{subsection3.2}

\begin{figure}[htbp]
 \centering\includegraphics[width=15cm,clip]{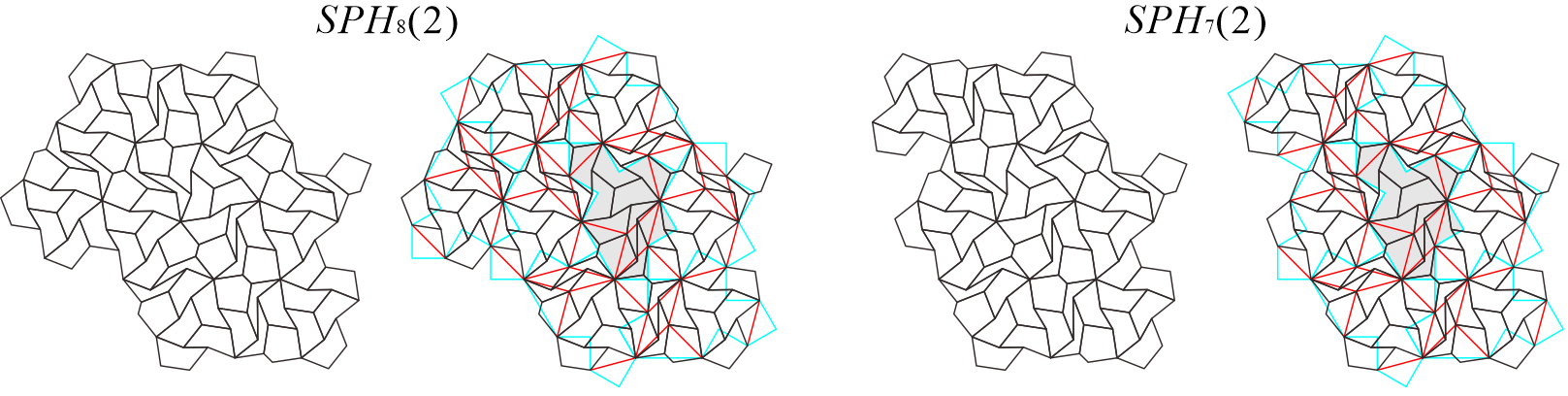} 
  \caption{{\small 
$SPH_{8}(2)$ and $SPH_{7}(2)$.
} 
\label{Fig.3-5}
}
\end{figure}

In the second pattern of $T_{h}$ (the second pattern when $T_{h}$ with \mbox{Tile$(1, 1)$} is 
converted into a tiling with three types of pentagons) in the case where the original rhombus 
is not subdivided into smaller similar rhombuses, all pentagons in the pairs corresponding to 
rhombuses with an acute angle of $30^ \circ$ are marked without asterisks (the anterior side in 
Figure~\ref{Fig.2-2}). The clusters formed by the three types of pentagons based on clusters 
$H_{8}$ and $H_{7}$ (see Figure~\ref{Fig.3-1}) for this pattern are shown in Figure~\ref{Fig.3-5}. 
Let these clusters be $SPH_{8}(2)$ and $SPH_{7}(2)$.

\begin{figure}[htbp]
 \centering\includegraphics[width=15cm,clip]{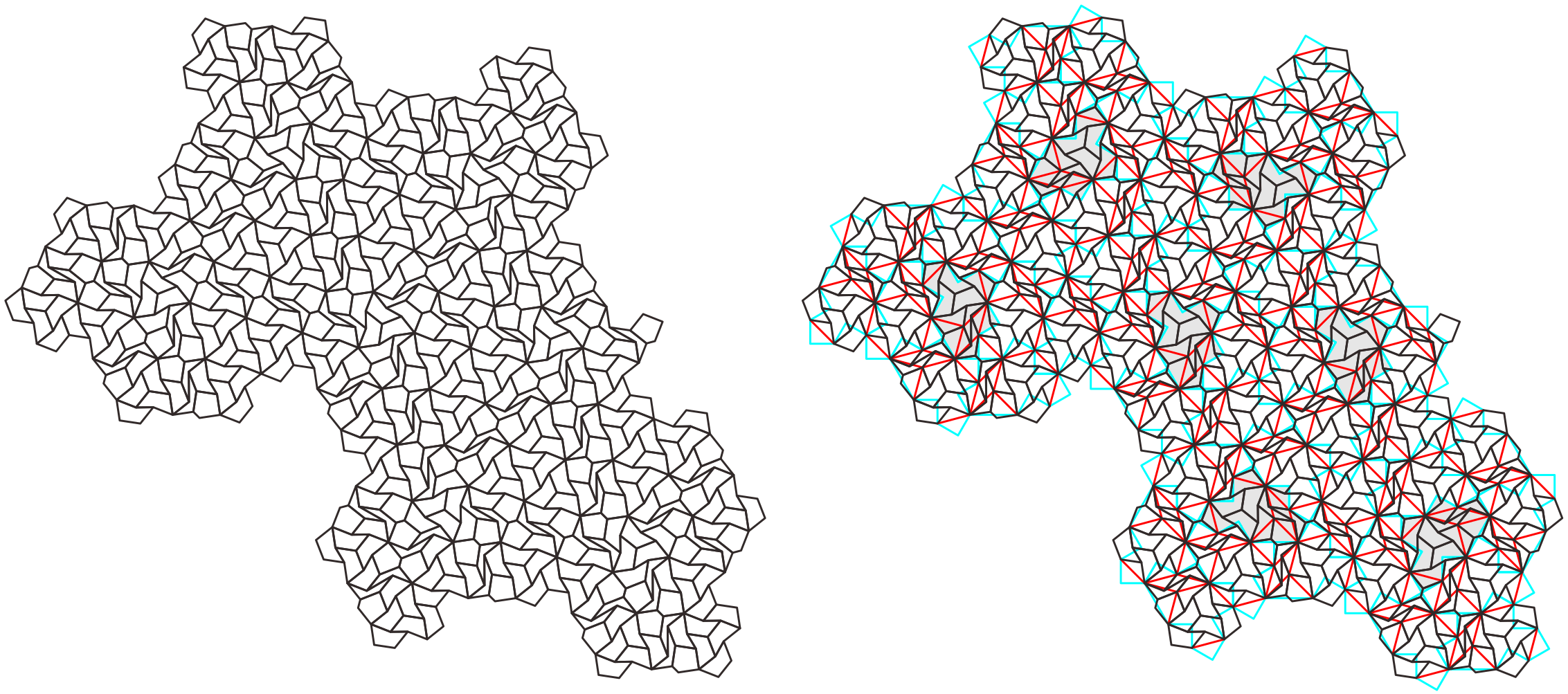} 
  \caption{{\small 
$SPH_{8}(3)$.
} 
\label{Fig.3-6}
}
\end{figure}

\begin{figure}[htbp]
 \centering\includegraphics[width=15cm,clip]{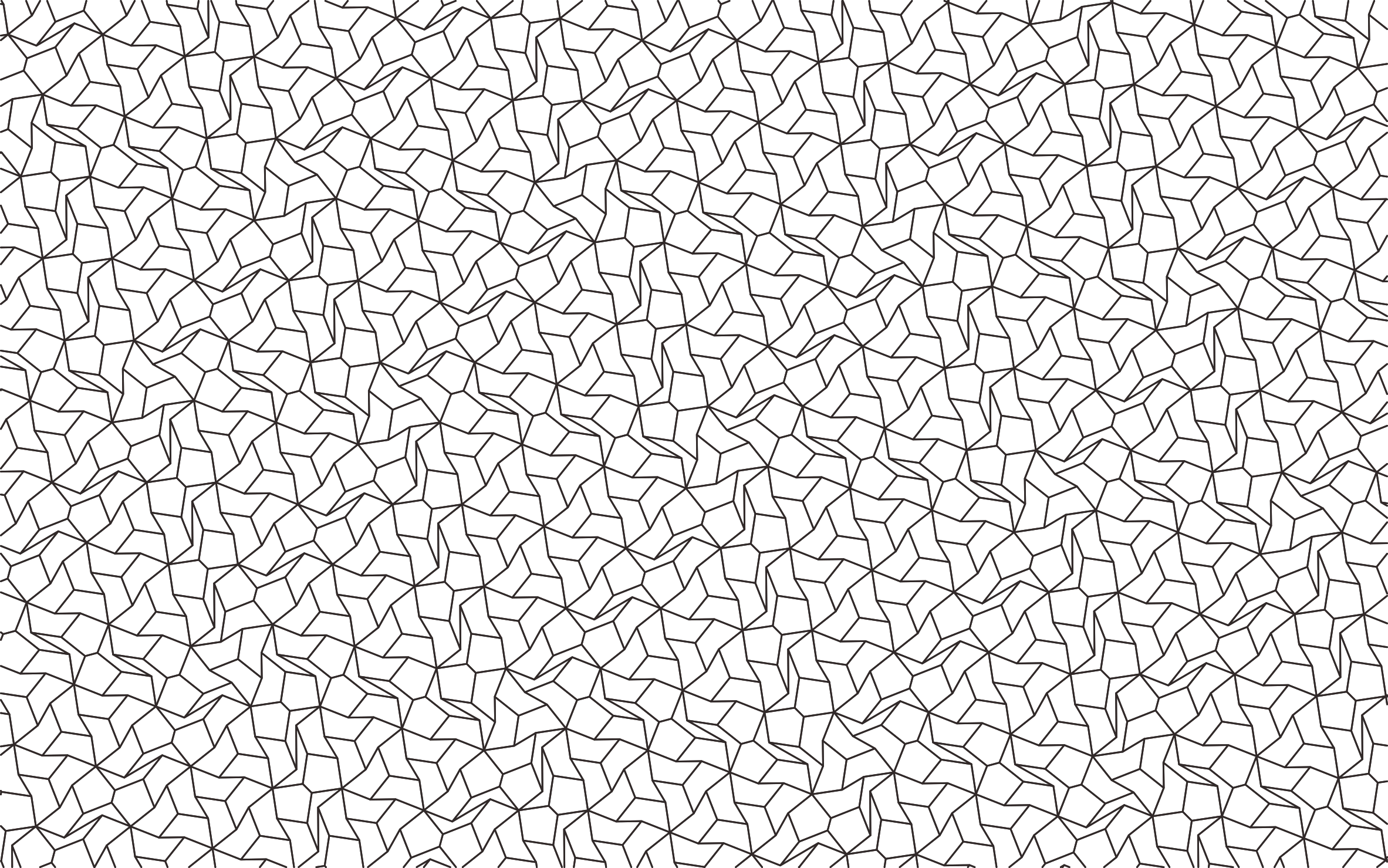} 
  \caption{{\small 
Second pattern of $T_{h}$ with three types of pentagons in the case where the original 
rhombus is not subdivided into smaller similar rhombuses.
} 
\label{Fig.3-7}
}
\end{figure}

Figure~\ref{Fig.3-6} shows cluster $SPH_{8}(3)$ of the next substitution step formed by 
$SPH_{8}(2)$ and $SPH_{7}(2)$. Figure~\ref{Fig.3-7} shows the second pattern of $T_{h}$ with 
three types of pentagons in the case where the original rhombus is not subdivided into smaller 
similar rhombuses, which was generated using $SPH_{8}(2)$ and $SPH_{7}(2)$ in conjunction 
with the substitution method described in \cite{Smith_2024a}.

\vspace{10mm}

\renewcommand{\figurename}{{\small Figure}}

\renewcommand{\figurename}{{\small Figure}}

\renewcommand{\figurename}{{\small Figure}}


\section{Conversion of non-periodic tiling $T_{s}$ in the case where the 
original rhombus is not subdivided into smaller similar rhombuses}
\label{section4}

\begin{figure}[htbp]
 \centering\includegraphics[width=15cm,clip]{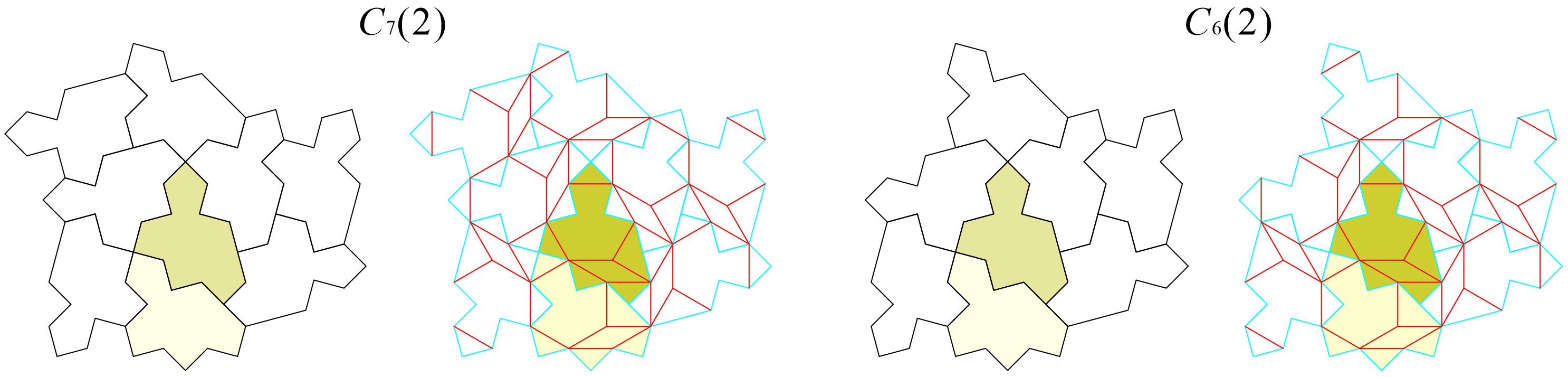} 
  \caption{{\small 
Clusters in Step 2 of Figure A.1 in \cite{Smith_2024b}.
} 
\label{Fig.4-1}
}
\end{figure}

When $T_{s}$ with \mbox{Tile$(1, 1)$} is converted into a tiling with three types of pentagons 
without subdivision of the original rhombus into smaller similar rhombuses (i.e., in the case where 
the original rhombus is not subdivided into smaller similar rhombuses), two patterns are observed: 
(i) the pentagons in pairs corresponding to rhombuses with an acute angle of $30^ \circ$ in 
certain locations are marked without asterisks (the anterior side in Figure~\ref{Fig.2-2}), 
whereas the majority of the pentagons in pairs corresponding to rhombuses with an acute angle 
of $30^ \circ$ in all other locations are marked with asterisks (the posterior side in 
Figure~\ref{Fig.2-2}); and (ii) the opposite configuration. Additionally, the patterns on 
\mbox{Tile$(1, 1)$} for forming each pattern of the pentagonal tilings can be classified into 
multiple types. In this study, we do not classify \mbox{Tile$(1, 1)$} based on different patterns, 
nor do we use \mbox{Tile$(1, 1)$} to generate tilings. In this study, the conversion results were 
obtained by generating tilings using the cluster in Step 2a of Figure A.1 in \cite{Smith_2024b} (see 
Figure~\ref{Fig.4-1}, denoted as $C_{7}(2)$) and the cluster obtained by removing one 
\mbox{Tile$(1, 1)$} from $C_{7}(2)$ (see Figure~\ref{Fig.4-1}, denoted as $C_{6}(2)$), in conjunction 
with the substitution method described in \cite{Smith_2024b}. In other words, we generated 
non-periodic tilings with three types of pentagons using the clusters composed of three types 
of pentagons based on $C_{7}(2)$ and $C_{6}(2)$ formed by \mbox{Tile$(1, 1)$}, in conjunction 
with the substitution method for $C_{7}(2)$ and $C_{6}(2)$ described in \cite{Smith_2024b}.

\begin{figure}[htbp]
 \centering\includegraphics[width=15cm,clip]{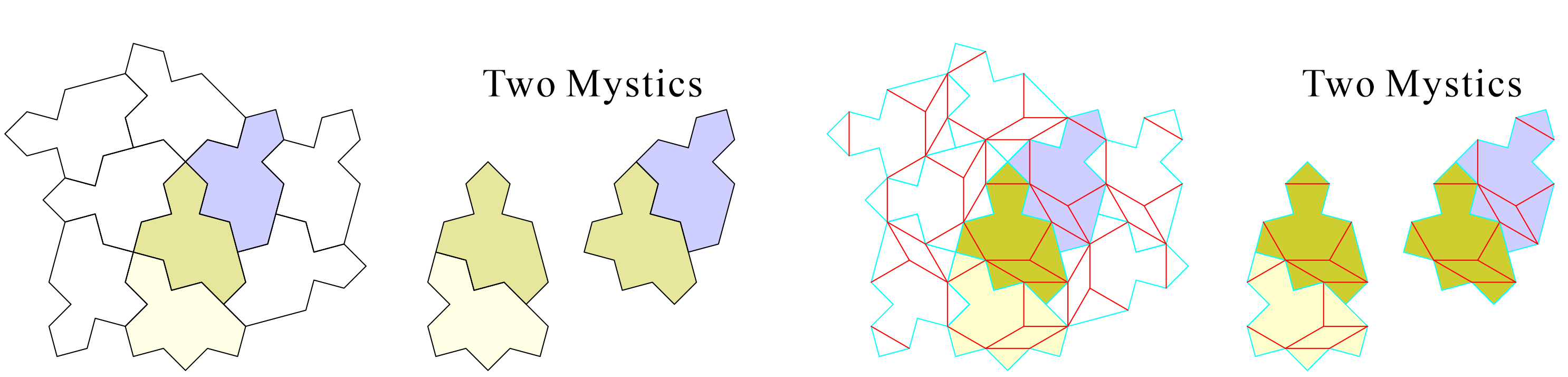} 
  \caption{{\small 
Two Mystic regions that formed by three \mbox{Tile$(1, 1)$} in $C_{7}(2)$.
} 
\label{Fig.4-2}
}
\end{figure}

The structure formed by the two colored \mbox{Tile$(1, 1)$} in Figure 4-1 is referred to as 
the ``Mystic'' region in \cite{Smith_2024b}. As shown in Figure~\ref{Fig.4-2}, we assumed that 
two locations correspond to the Mystic structure in $C_{7}(2)$ (i.e., two Mystic regions formed 
by three \mbox{Tile$(1, 1)$}).

Based on the substitution method described in the Appendix of \cite{Smith_2024b}, the tiles 
are switched between the anterior and posterior sides (i.e., the tiles are reversed) at each 
substitution step. For example, if $C_{7}(3)$ corresponding to the cluster in Step 3 is formed 
by performing a substitution using clusters $C_{7}(2)$ and $C_{6}(2)$ in Figure~\ref{Fig.4-1}, 
then all \mbox{Tile$(1, 1)$} in $C_{7}(3)$ will be at the posterior side in Figure~\ref{Fig.1-1}. 
However, in the figures presented herein, we did not adopt a notation to indicate that the 
tiles are reversed at each substitution step. Instead, we adjusted the figures so that 
\mbox{Tile$(1, 1)$} in the tilings and clusters were consistently oriented toward the anterior 
side shown in Figure~\ref{Fig.1-1}.

\renewcommand{\figurename}{{\small Figure}}

\renewcommand{\figurename}{{\small Figure}}


\subsection{First pattern of $T_{s}$ in the case where the original rhombus is not 
subdivided into smaller similar rhombuses}
\label{subsection4.1}

\begin{figure}[htbp]
 \centering\includegraphics[width=15cm,clip]{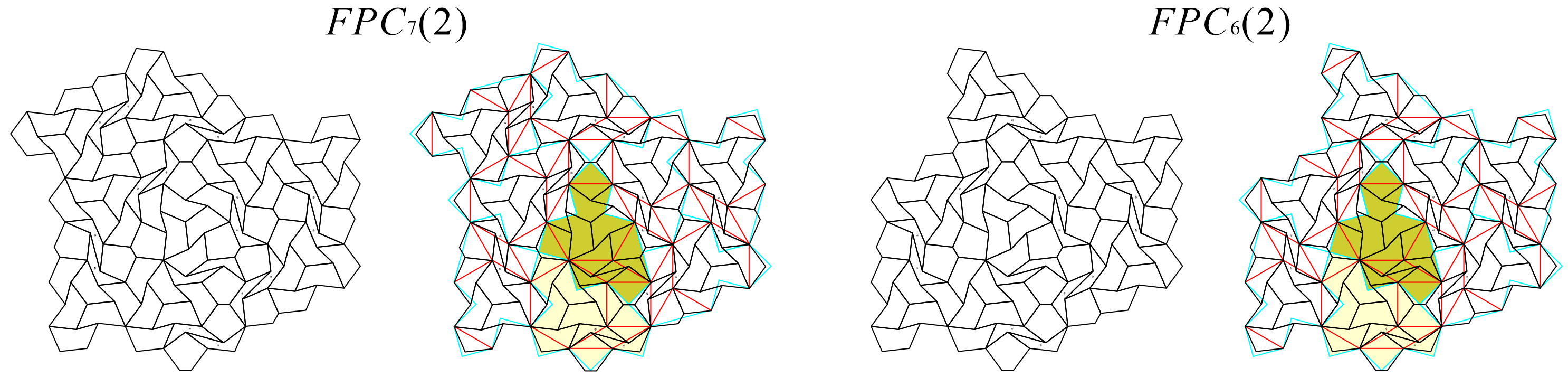} 
  \caption{{\small 
$FPC_{7}(2)$ and $FPC_{6}(2)$.
} 
\label{Fig.4-3}
}
\end{figure}

In the first pattern of $T_{s}$ (the first pattern when $T_{s}$ with \mbox{Tile$(1, 1)$} is converted 
into a tiling with three types of pentagons) in the case where the original rhombus is not subdivided 
into smaller similar rhombuses, clusters $FPC_{7}(2)$ and $FPC_{6}(2)$, consisting of three types 
of pentagons based on $C_{7}(2)$ and $C_{6}(2)$, respectively, are shown in Figure~\ref{Fig.4-3}.

\begin{figure}[htbp]
 \centering\includegraphics[width=15cm,clip]{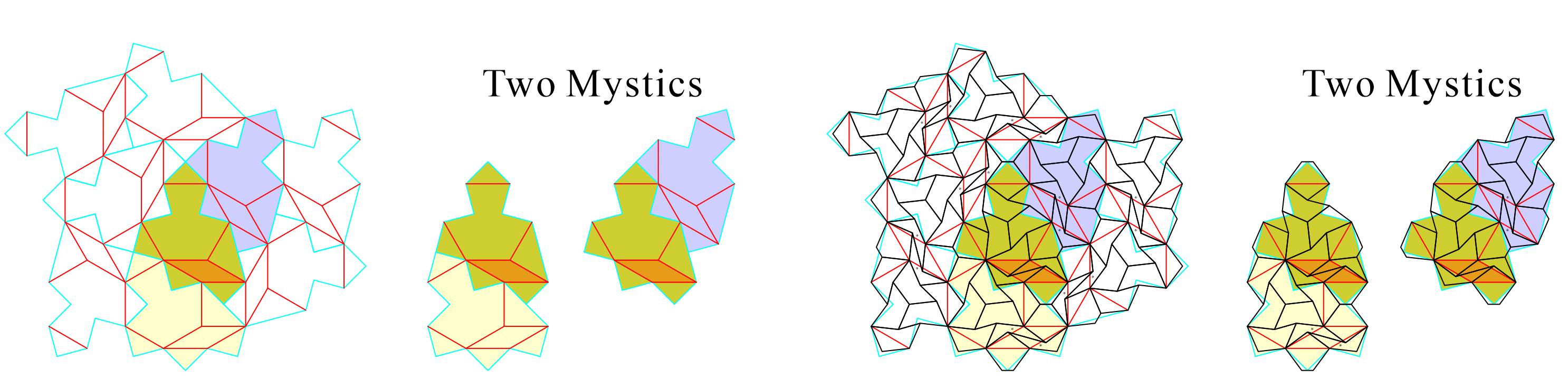} 
  \caption{{\small 
  Rhombus with an acute angle of $30^ \circ$ within \mbox{Tile$(1, 1)$}, which 
appears in the overlap of two Mystic regions in $C_{7}(2)$ and $FPC_{7}(2)$.
} 
\label{Fig.4-4}
}
\end{figure}

In the first pattern of $T_{s}$, the pentagons in the pair corresponding to the rhombus with 
an acute angle of $30^ \circ$ (see the rhombus with an orange filter in Figure~\ref{Fig.4-4}) 
within \mbox{Tile$(1, 1)$}\footnote{ 
The overlapping \mbox{Tile$(1, 1)$} corresponds to the ``odd tile'' in \cite{Smith_2024b}.
}, which appears in the overlap of two Mystic regions in $C_{7}(2)$ and $C_{6}(2)$, are marked 
without asterisks (the anterior side in Figure~\ref{Fig.2-2}). Meanwhile, the other pentagons in pairs 
corresponding to the remaining rhombuses with an acute angle of $30^ \circ$ are marked 
with asterisks (the posterior side in Figure~\ref{Fig.2-2})\footnote{ 
In this version (Ver. 3), the relationship between the presence and absence of asterisks, as 
used in the explanation of Ver. 1, is reversed based on the contents of Section~\ref{section3}.
\label{footnote08}
}.

\begin{figure}[htbp]
 \centering\includegraphics[width=15cm,clip]{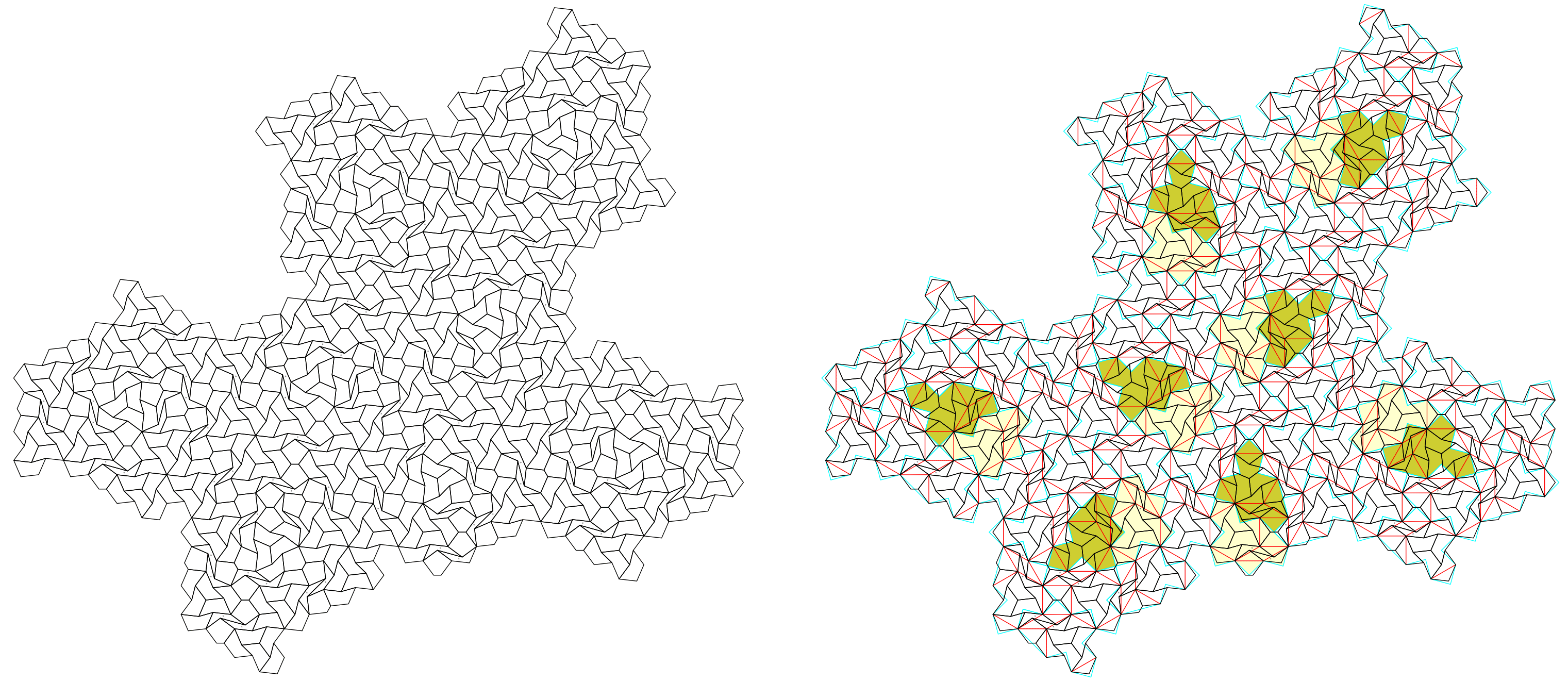} 
  \caption{{\small 
$FPC_{7}(3)$ corresponding to the cluster in Step 3a of Figure A.1 
in \cite{Smith_2024b} (Note: The clusters in this figure are not reversed).
} 
\label{Fig.4-5}
}
\end{figure}

\begin{figure}[htbp]
 \centering\includegraphics[width=15cm,clip]{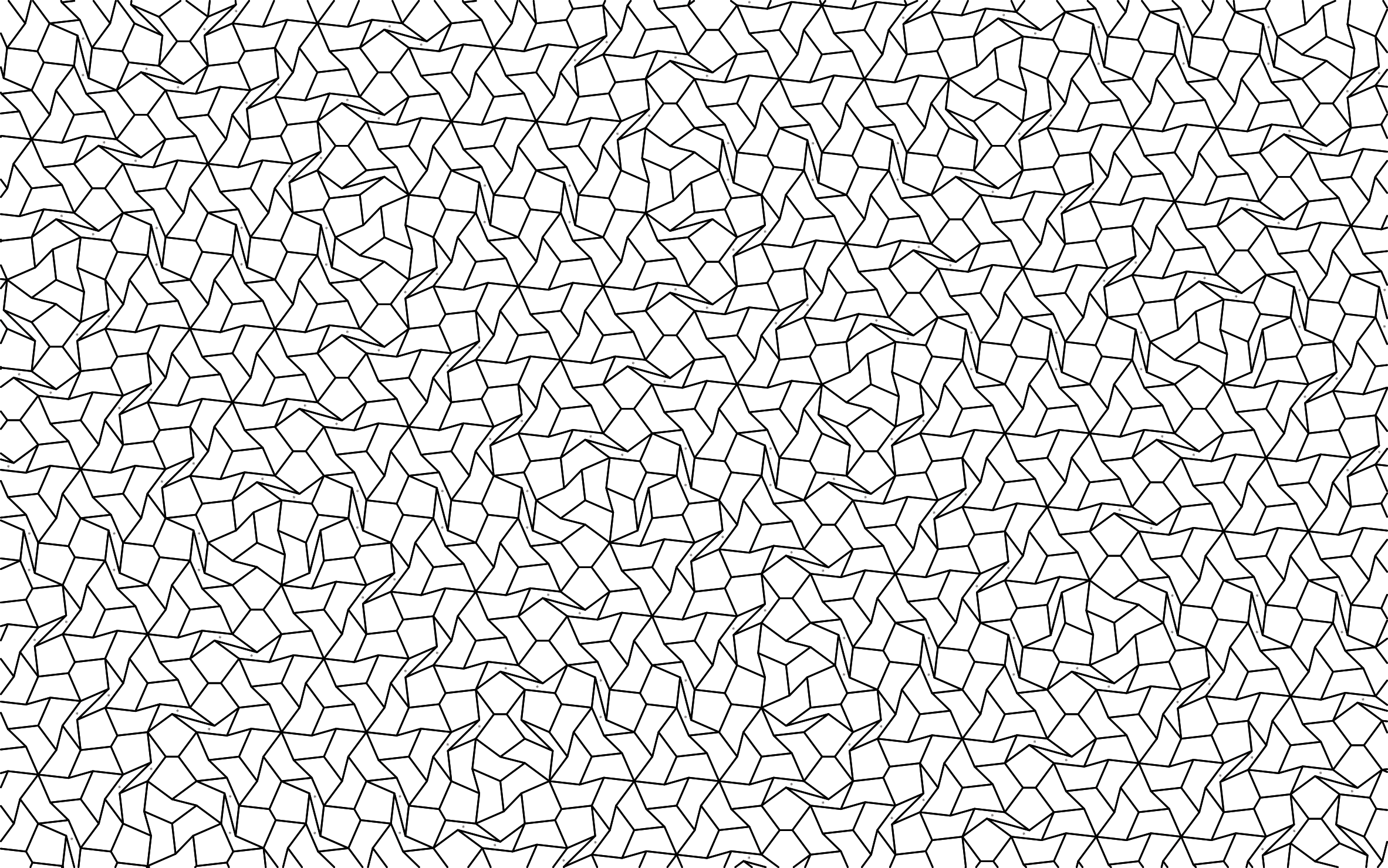} 
  \caption{{\small 
First pattern of $T_{s}$ with three types of pentagons in the case where the original 
rhombus is not subdivided into smaller similar rhombuses.
} 
\label{Fig.4-6}
}
\end{figure}

Figure~\ref{Fig.4-5} shows cluster $FPC_{7}(3)$, which corresponds to the cluster in Step 3a of 
Figure A.1 in \cite{Smith_2024b} and is formed by $FPC_{7}(2)$ and $FPC_{6}(2)$\footnote{ 
As mentioned above, in the figures presented herein, we did not adopt a notation to indicate 
that the tiles are reversed at each substitution step. Note that the clusters in Figure~\ref{Fig.4-5} 
and the cluster in Step 3a of Figure A.1 in \cite{Smith_2024b} are a reflective relationship 
(reversed direction). \label{footnote09}
}. Figure~\ref{Fig.4-6} shows the first pattern of $T_{s}$ with three types of pentagons in the 
case where the original rhombus is not subdivided into smaller similar rhombuses, which was 
generated using $FPC_{7}(2)$ and $FPC_{6}(2)$ in conjunction with the substitution method 
described in \cite{Smith_2024b}.

\renewcommand{\figurename}{{\small Figure}}

\renewcommand{\figurename}{{\small Figure}}

\renewcommand{\figurename}{{\small Figure}}

\renewcommand{\figurename}{{\small Figure}}


\subsection{Second pattern of $T_{s}$ in the case where the original rhombus is not 
subdivided into smaller similar rhombuses}
\label{subsection4.2}

\begin{figure}[htbp]
 \centering\includegraphics[width=15cm,clip]{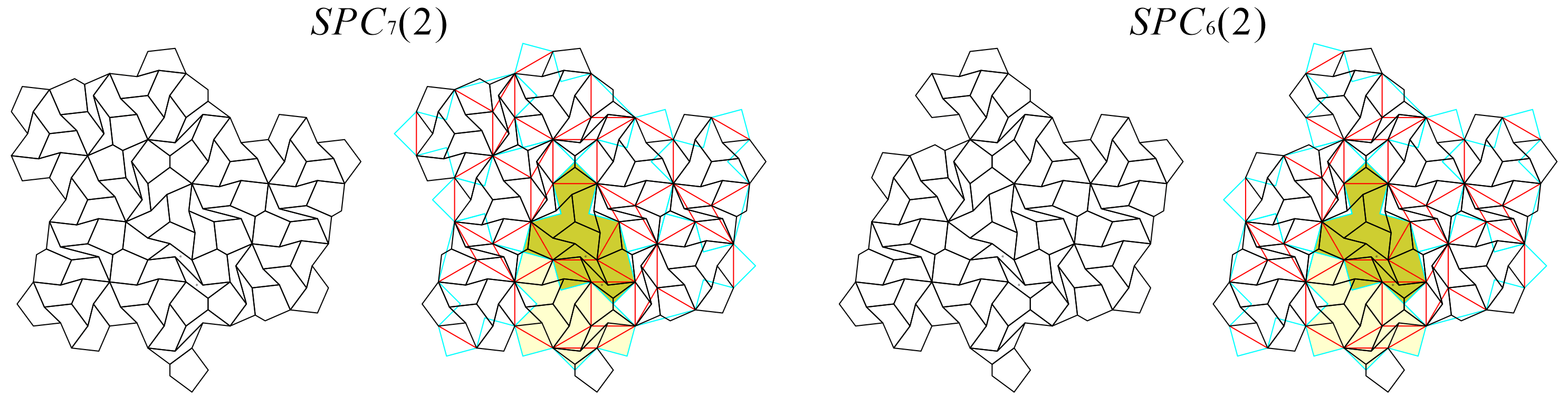} 
  \caption{{\small 
$SPC_{7}(2)$ and $SPC_{6}(2)$.
} 
\label{Fig.4-7}
}
\end{figure}

In the second pattern of $T_{s}$ (the second pattern when $T_{s}$ with \mbox{Tile$(1, 1)$} is 
converted into a tiling with three types of pentagons) in the case where the original rhombus 
is not subdivided into smaller similar rhombuses, clusters $SPC_{7}(2)$ and $SPC_{6}(2)$, 
consisting of three types of pentagons based on $C_{7}(2)$ and $C_{6}(2)$, respectively, 
are shown in Figure~\ref{Fig.4-7}.

In the second pattern of $T_{s}$, the pentagons in the pair corresponding to the rhombus 
with an acute angle of $30^ \circ$ (see the rhombus with an orange filter in Figure~\ref{Fig.4-4}) 
within \mbox{Tile$(1, 1)$}, which appears in the overlap of two Mystic regions in $C_{7}(2)$ and 
$C_{6}(2)$, are marked with asterisks (the posterior side in Figure~\ref{Fig.2-2}). Meanwhile, the 
other pentagons in pairs corresponding to the remaining rhombuses with an acute angle of 
$30^ \circ$ are marked without asterisks (the anterior side in Figure~\ref{Fig.2-2})\footref{footnote08}.

\begin{figure}[htbp]
 \centering\includegraphics[width=15cm,clip]{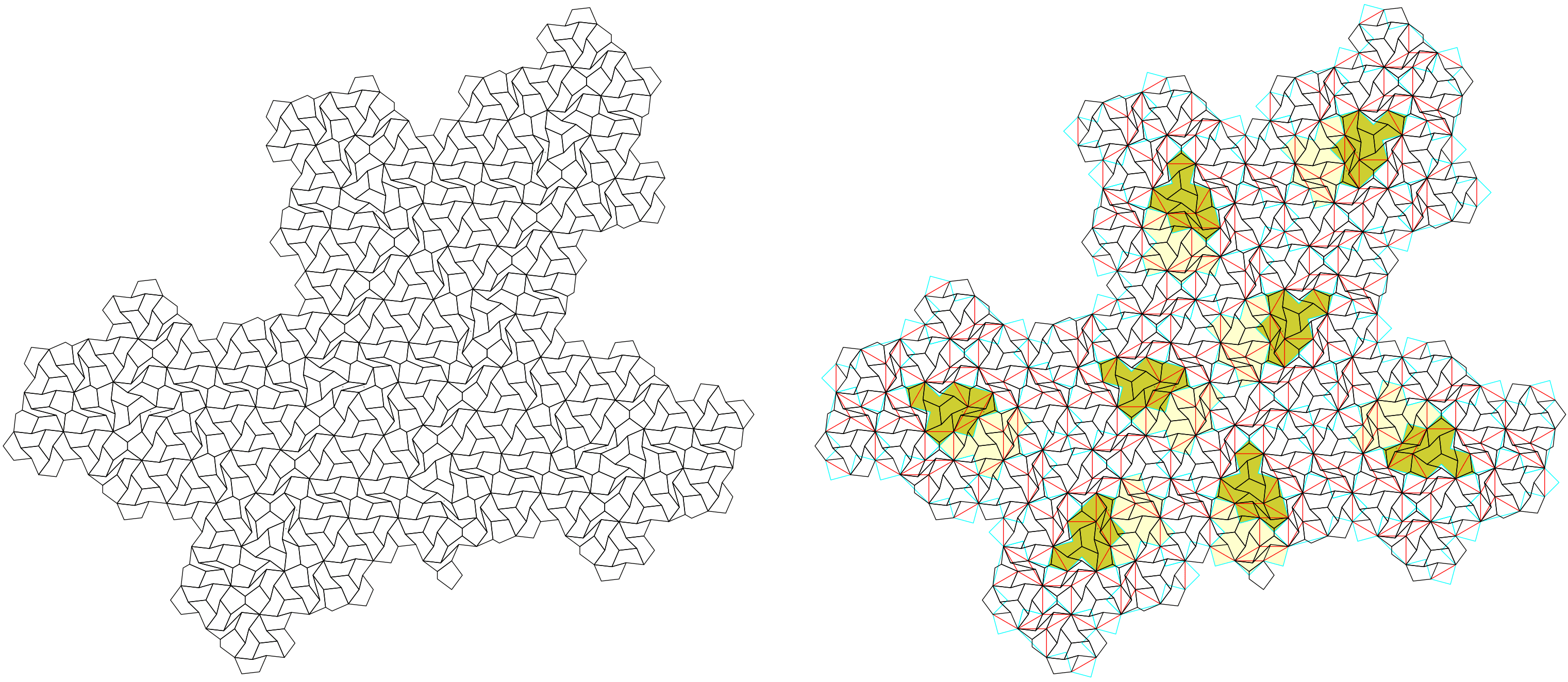} 
  \caption{{\small 
$SPC_{7}(3)$ corresponding to the cluster in Step 3a of Figure A.1 
in \cite{Smith_2024b} (Note: The clusters in this figure are not reversed).
} 
\label{Fig.4-8}
}
\end{figure}

\begin{figure}[htbp]
 \centering\includegraphics[width=15cm,clip]{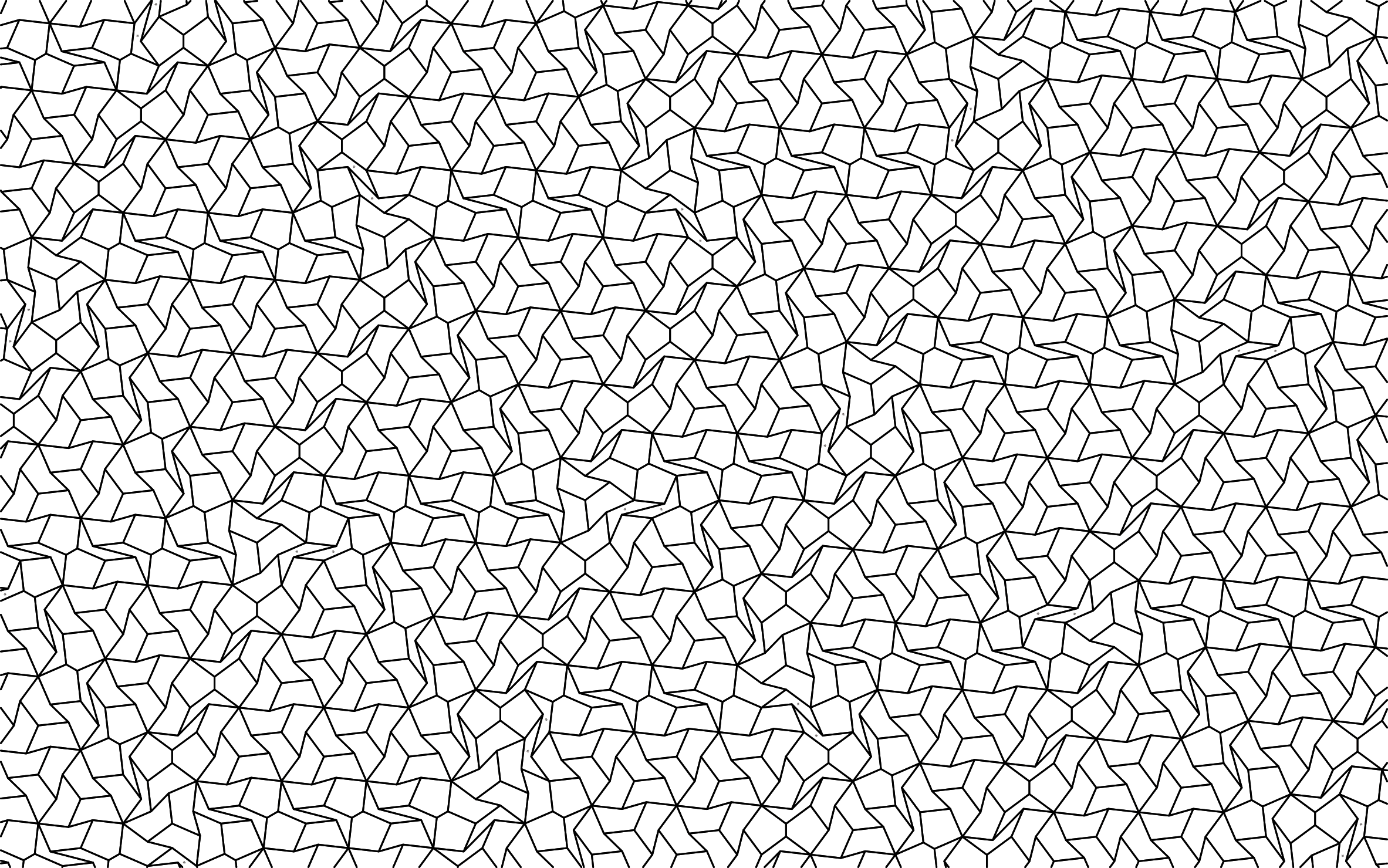} 
  \caption{{\small 
Second pattern of $T_{s}$ with three types of pentagons in the case where the 
original rhombus is not subdivided into smaller similar rhombuses.
} 
\label{Fig.4-9}
}
\end{figure}

Figure~\ref{Fig.4-8} shows cluster $SPC_{7}(3)$, which corresponds to the cluster in Step 3a 
of Figure A.1 in \cite{Smith_2024b} and is formed by $SPC_{7}(2)$ and $SPC_{6}(2)$\footref{footnote09}. 
Figure~\ref{Fig.4-9} shows the second pattern of $T_{s}$ with three types of pentagons in 
the case where the original rhombus is not subdivided into smaller similar rhombuses, 
which was generated using $SPC_{7}(2)$ and $SPC_{6}(2)$ in conjunction with the 
substitution method described in \cite{Smith_2024b}.

\renewcommand{\figurename}{{\small Figure}}

\renewcommand{\figurename}{{\small Figure}}

\renewcommand{\figurename}{{\small Figure}}


\subsection{Countless patterns of non-periodic tilings that \mbox{Tile$(1, 1)$} can 
generate and their conversion}
\label{subsection4.3}

\begin{figure}[t]
 \centering\includegraphics[width=15cm,clip]{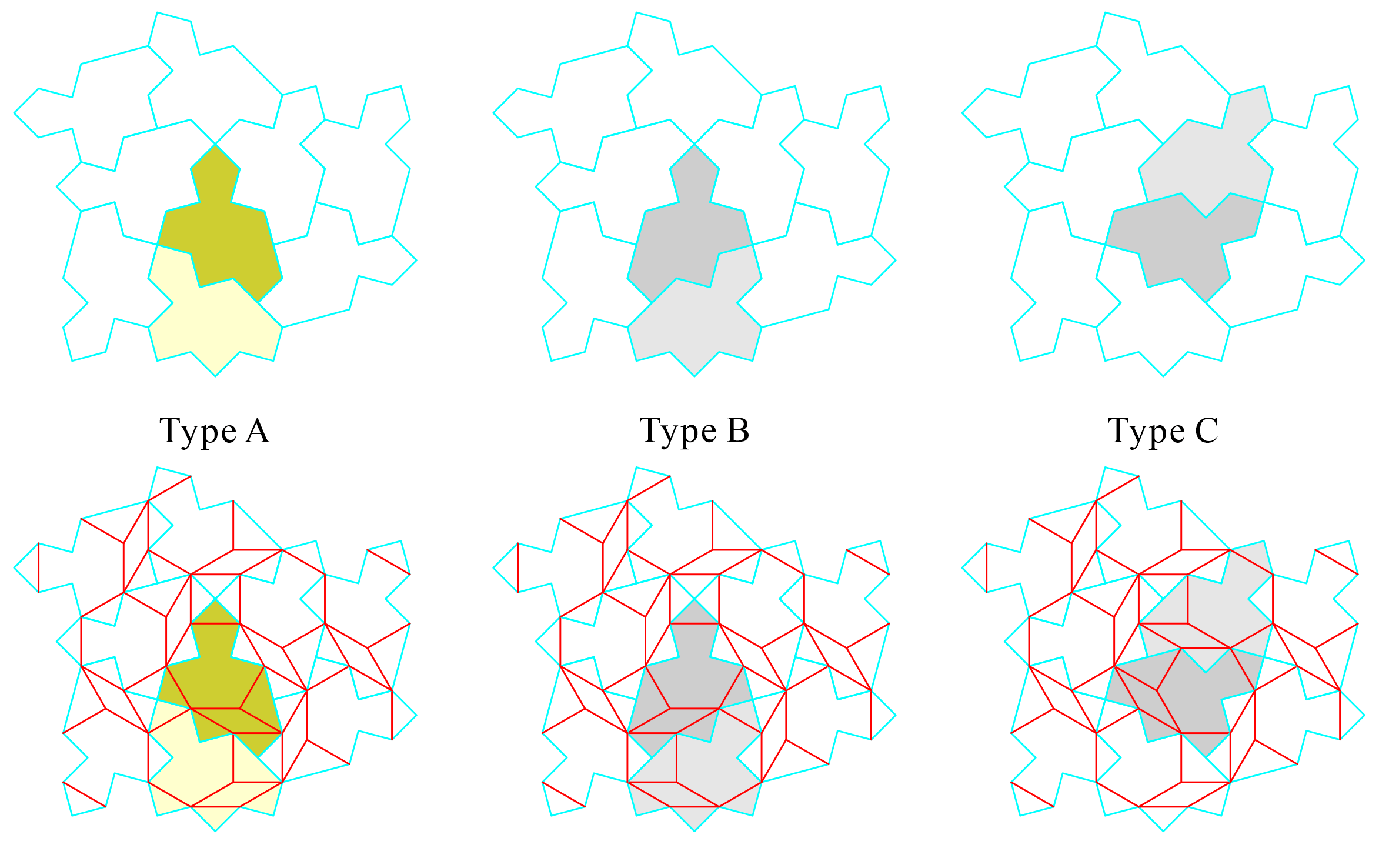} 
  \caption{{\small 
$C_{7}(2)$ and three types of patterns.
} 
\label{Fig.4-10}
}
\end{figure}

As shown in Figure~\ref{Fig.4-4}, the shape of the Mystic exhibits line symmetry; thus, 
its outline remains unchanged under reflection. Therefore, clusters $C_{7}(2)$ (and $C_{6}(2)$) 
can create three types of patterns, as shown in Figure~\ref{Fig.4-10}, without changing the 
outline. (Type A in Figure~\ref{Fig.4-10} corresponds to $C_{7}(2)$ of $T_{s}$ in Figure~\ref{Fig.4-1}. 
Type B and Type C in Figure~\ref{Fig.4-10} are clusters created by replacing the Mystic 
regions inside $C_{7}(2)$ with the Mystic formed by the reflected \mbox{Tile$(1, 1)$}, 
and their outlines are identical to that of Type A). We focused on the steps following Step 2a 
of the substitution method described in the Appendix of \cite{Smith_2024b} and noticed 
that non-periodic tilings can be generated using two types of tiles with the outlines of 
$C_{7}(2)$ and $C_{6}(2)$ in conjunction with the substitution method. This suggests that 
non-periodic tilings with clusters of Type B or Type C (as shown in Figure~\ref{Fig.4-10}) 
can be generated using the substitution method that generates $T_{s}$ using $C_{7}(2)$ 
and $C_{6}(2)$. Furthermore, by replacing the locations corresponding to the structure of 
Mystic in $T_{s}$, as generated using $C_{7}(2)$ and $C_{6}(2)$, with the Mystic formed 
by the reflected \mbox{Tile$(1, 1)$} (this replacement can be performed arbitrarily), 
we can generate countless patterns (design patterns created by the arrangement of 
polygonal tiles) of non-periodic tilings based on $T_{s}$, although they will be different 
from the original patterns of $T_{s}$. In other words, if the use of reflected tiles is allowed 
during tiling generation, then \mbox{Tile$(1, 1)$} can generate countless patterns of tilings 
(notably, \mbox{Tile$(1, 1)$} can generate countless patterns of non-periodic tilings distinct 
from $T_{h}$ and $T_{s}$).

\begin{figure}[htbp]
 \centering\includegraphics[width=15cm,clip]{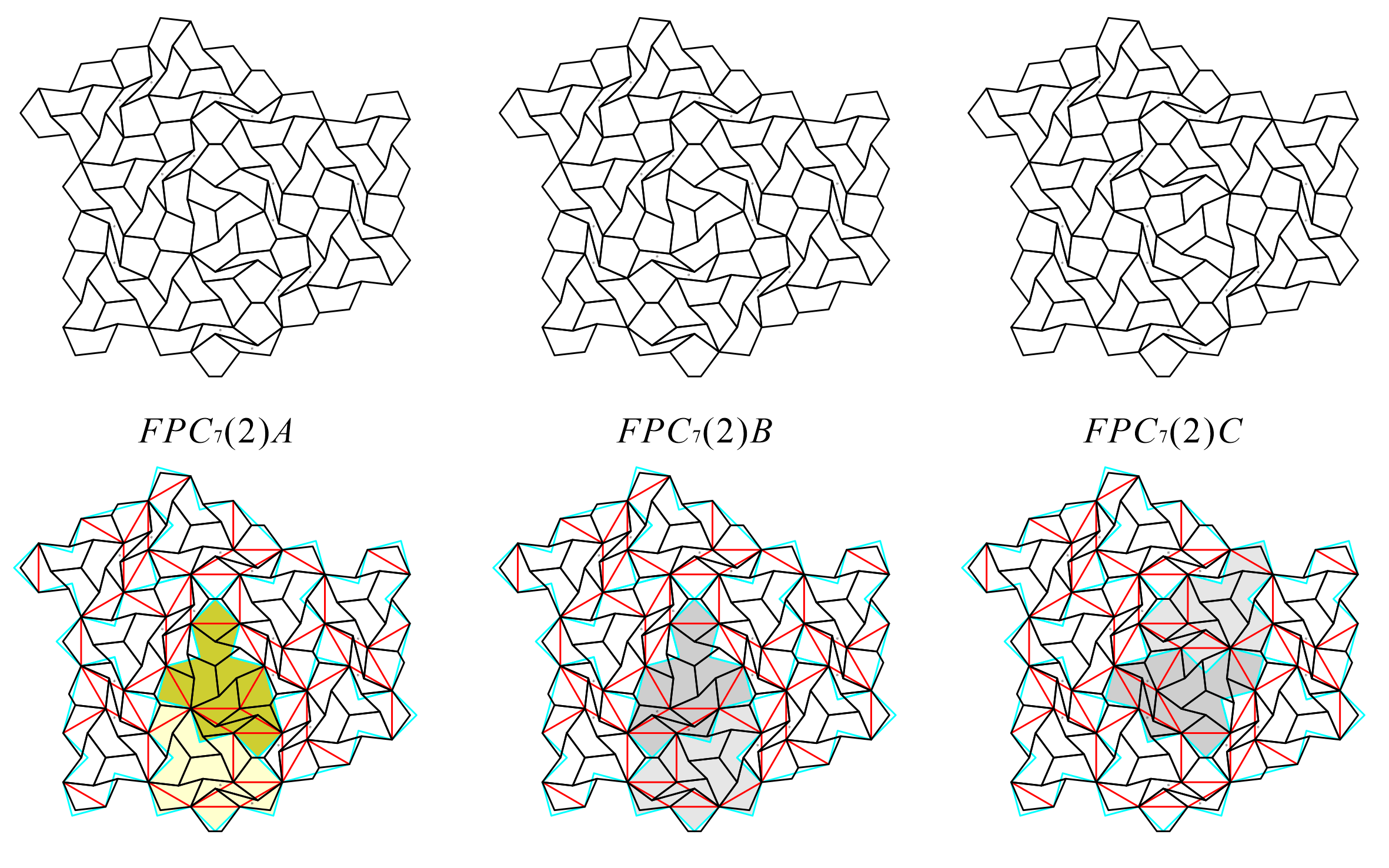} 
  \caption{{\small 
$FPC_{7}(2)A$, $FPC_{7}(2)B$, and $FPC_{7}(2)C$.
} 
\label{Fig.4-11}
}
\end{figure}

\begin{figure}[htbp]
 \centering\includegraphics[width=15cm,clip]{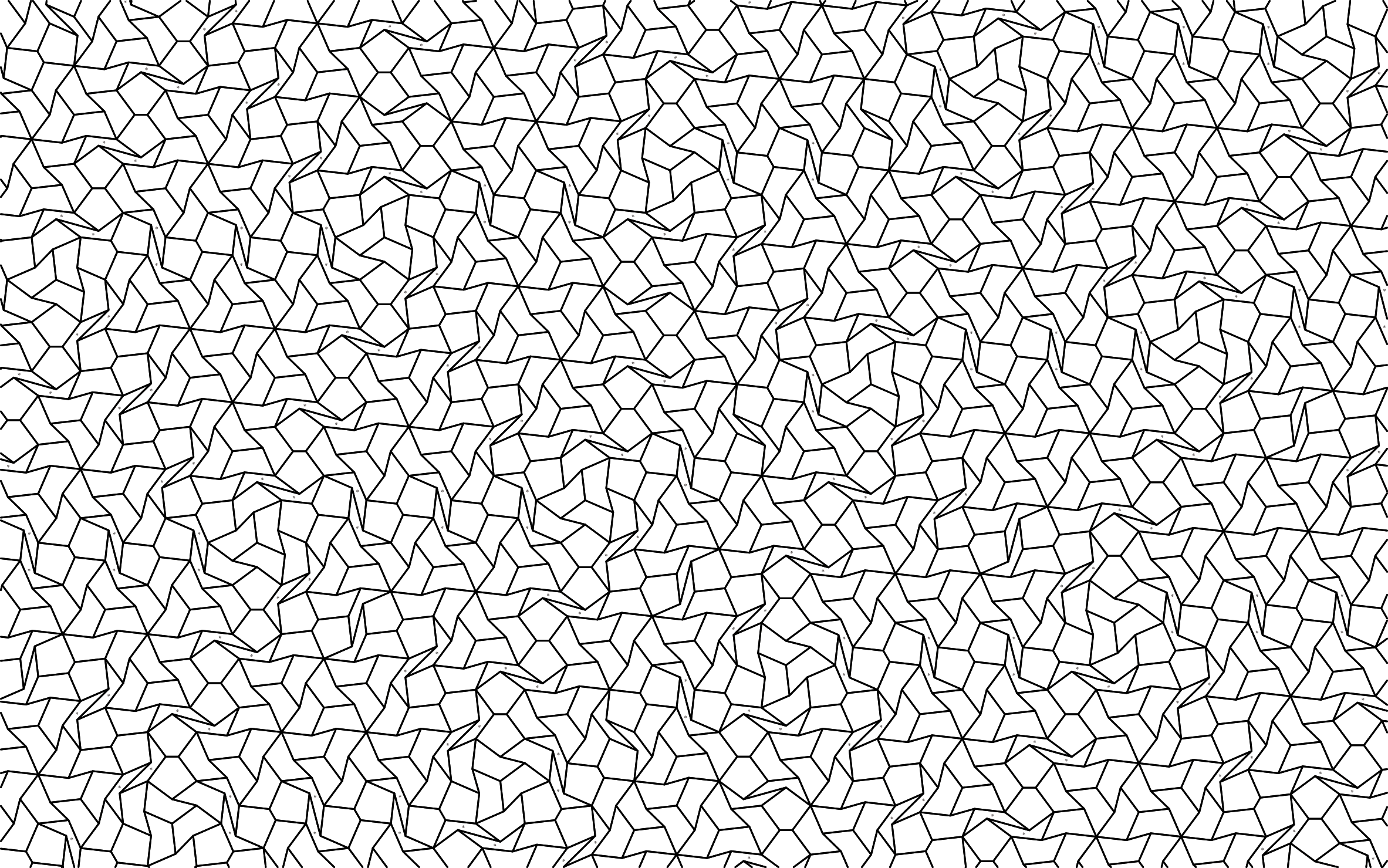} 
  \caption{{\small 
Non-periodic tiling with three types of pentagons generated 
using $FPC_{7}(2)B$ in conjunction with the substitution method described in \cite{Smith_2024b}.
} 
\label{Fig.4-12}
}
\end{figure}

\begin{figure}[htbp]
 \centering\includegraphics[width=15cm,clip]{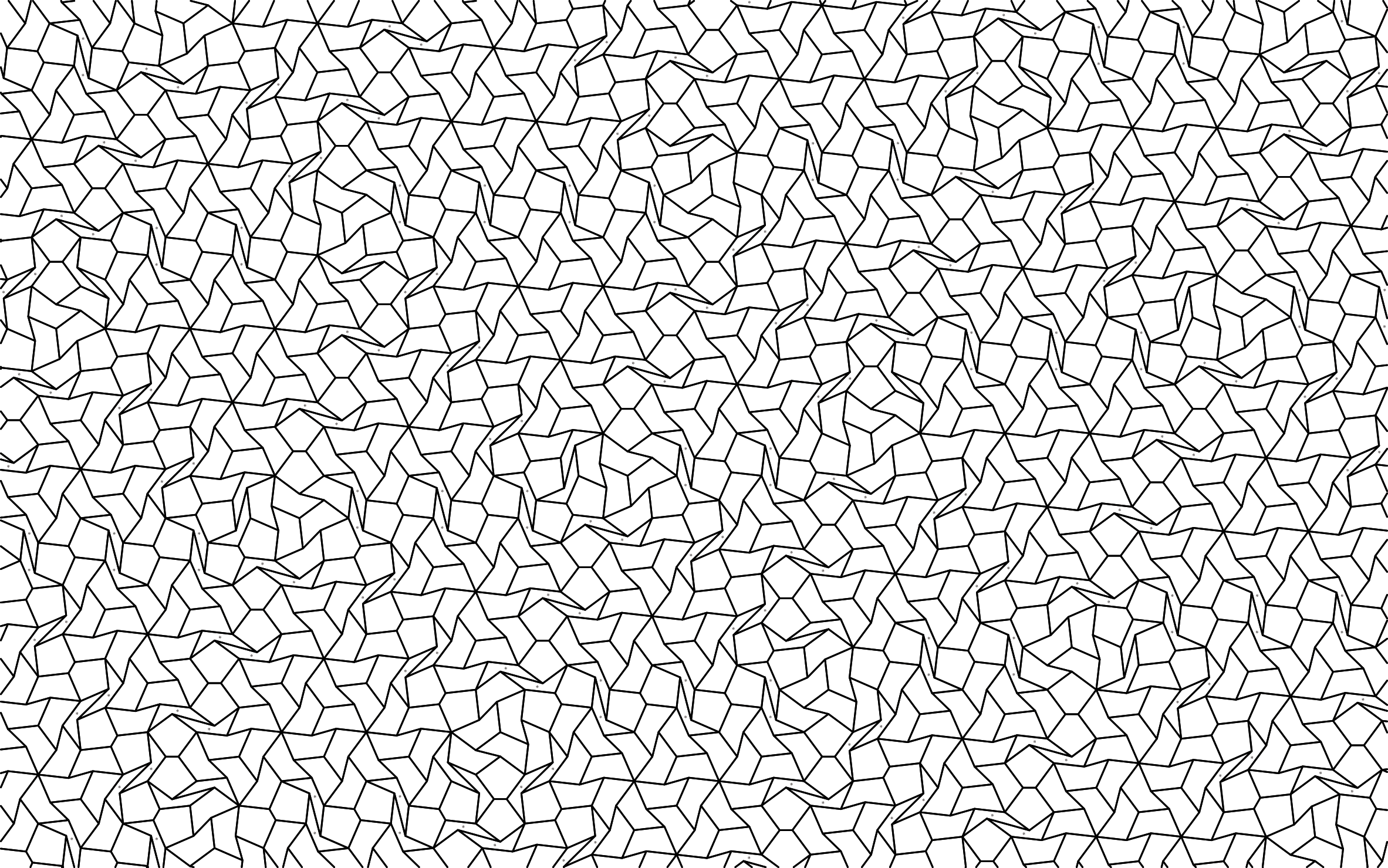} 
  \caption{{\small 
Non-periodic tiling with three types of pentagons generated 
using $FPC_{7}(2)C$ in conjunction with the substitution method described in \cite{Smith_2024b}.
} 
\label{Fig.4-13}
}
\end{figure}

Figure~\ref{Fig.4-11} shows clusters $FPC_{7}(2)A$, $FPC_{7}(2)B$, and $FPC_{7}(2)C$, 
which are based on the first pattern corresponding to Types A, B, and C in Figure~\ref{Fig.4-10}. 
(Note: $FPC_{7}(2)A$ corresponds to $FPC_{7}(2)$ in Figure~\ref{Fig.4-3}. The outlines of 
$FPC_{7}(2)B$ and $FPC_{7}(2)C$ are identical to that of $FPC_{7}(2)A$). Figure~\ref{Fig.4-12}  
shows a non-periodic tiling with three different pentagons, which was generated using $FPC_{7}(2)B$ 
in conjunction with the substitution method described in \cite{Smith_2024b}. Figure~\ref{Fig.4-13} 
shows a non-periodic tiling with three types of pentagons, which was generated using $FPC_{7}(2)C$ 
in conjunction with the substitution method described in \cite{Smith_2024b}.

\renewcommand{\figurename}{{\small Figure}}

\renewcommand{\figurename}{{\small Figure}}

\renewcommand{\figurename}{{\small Figure}}

\renewcommand{\figurename}{{\small Figure}}

\begin{figure}[htbp]
 \centering\includegraphics[width=15cm,clip]{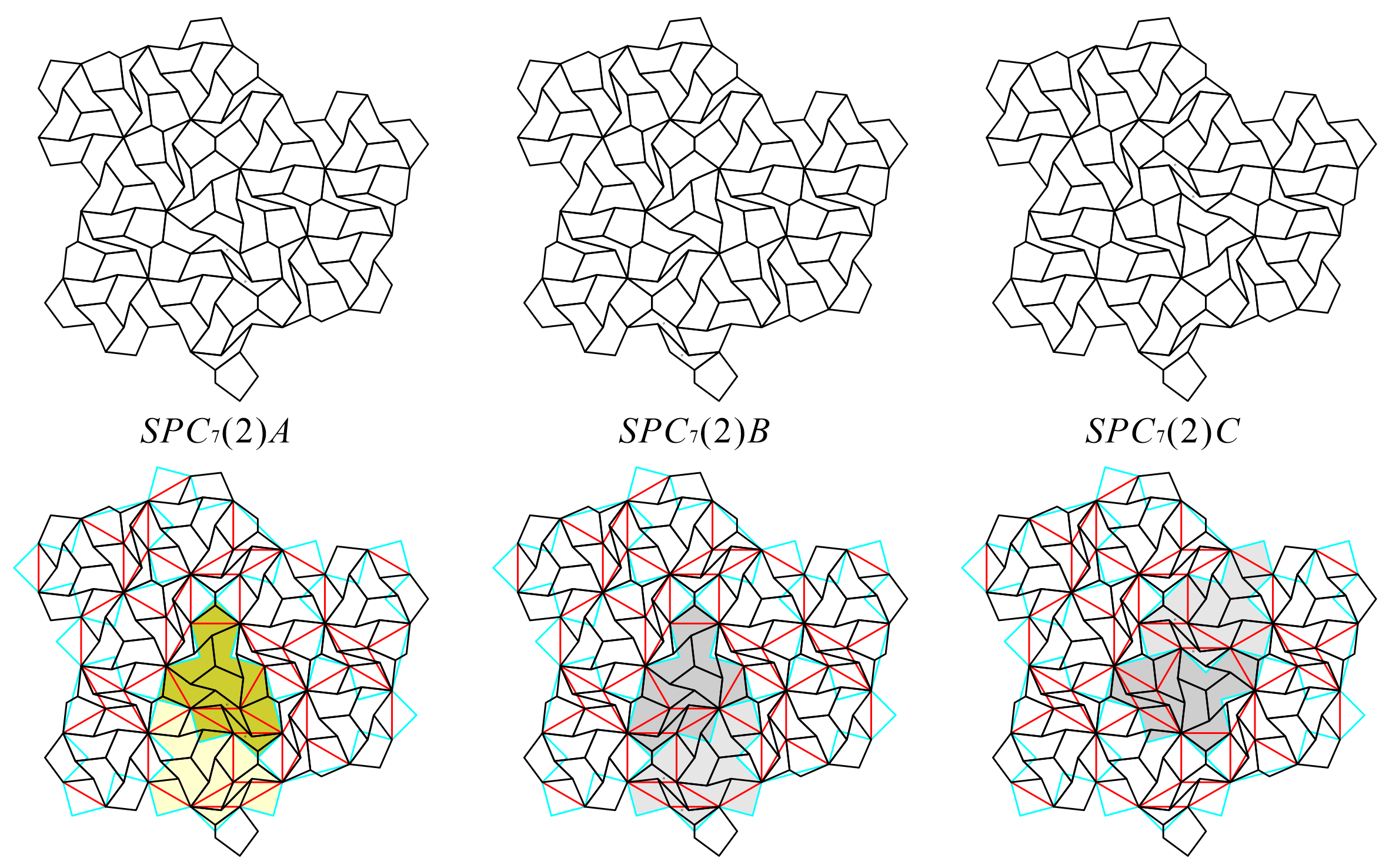} 
  \caption{{\small 
$SPC_{7}(2)A$, $SPC_{7}(2)B$, and $SPC_{7}(2)C$.
} 
\label{Fig.4-14}
}
\end{figure}

Figure~\ref{Fig.4-14} shows clusters $SPC_{7}(2)A$, $SPC_{7}(2)B$, and $SPC_{7}(2)C$, which are 
based on the second pattern corresponding to Types A, B, and C in Figure~\ref{Fig.4-10}. 
(Note:  $SPC_{7}(2)A$ corresponds to  $SPC_{7}(2)$, in Figure~\ref{Fig.4-7}. The outlines of 
 $SPC_{7}(2)B$ and  $SPC_{7}(2)C$ are identical to that of $SPC_{7}(2)A$). Similar to the first 
 pattern, by using $SPC_{7}(2)B$ or $SPC_{7}(2)C$ in conjunction with the substitution method 
 described in \cite{Smith_2024b}, non-periodic tilings can be generated (the illustrations of these 
 tilings are available on the author's website \cite{Sugimoto_site_Ts}).

\renewcommand{\figurename}{{\small Figure}}

Hence, we consider that countless patterns of non-periodic tilings with three types of pentagons 
can be generated based on $T_{s}$ (although the patterns fall into the first and second series). 
Based on Figure~\ref{Fig.4-11}, we assume that the first pattern series of $T_{s}$, where the 
original rhombus is not subdivided into smaller similar rhombuses, represents the case in which 
the pentagons at one of the three locations corresponding to rhombuses with an acute angle 
of $30^ \circ$, which are contained in three \mbox{Tile$(1, 1)$} used in the two Mystic regions 
in $C_{7}(2)$ (and $C_{6}(2)$), are marked without asterisks, whereas the pentagons at locations 
corresponding to the other rhombuses with an acute angle of $30^ \circ$ are marked with 
asterisks. Based on Figure~\ref{Fig.4-14}, we assume that the second pattern series of $T_{s}$, 
where the original rhombus is not subdivided into smaller similar rhombuses, represents the 
case in which the pentagons at one of the three locations corresponding to rhombuses with 
an acute angle of $30^ \circ$, which are contained in three \mbox{Tile$(1, 1)$} used in the 
two Mystic regions in $C_{7}(2)$ (and $C_{6}(2)$), are marked with asterisks, whereas the 
pentagons at locations corresponding to the other rhombuses with an acute angle of $30^ \circ$ 
are marked without asterisks.

\begin{figure}[htbp]
 \centering\includegraphics[width=15cm,clip]{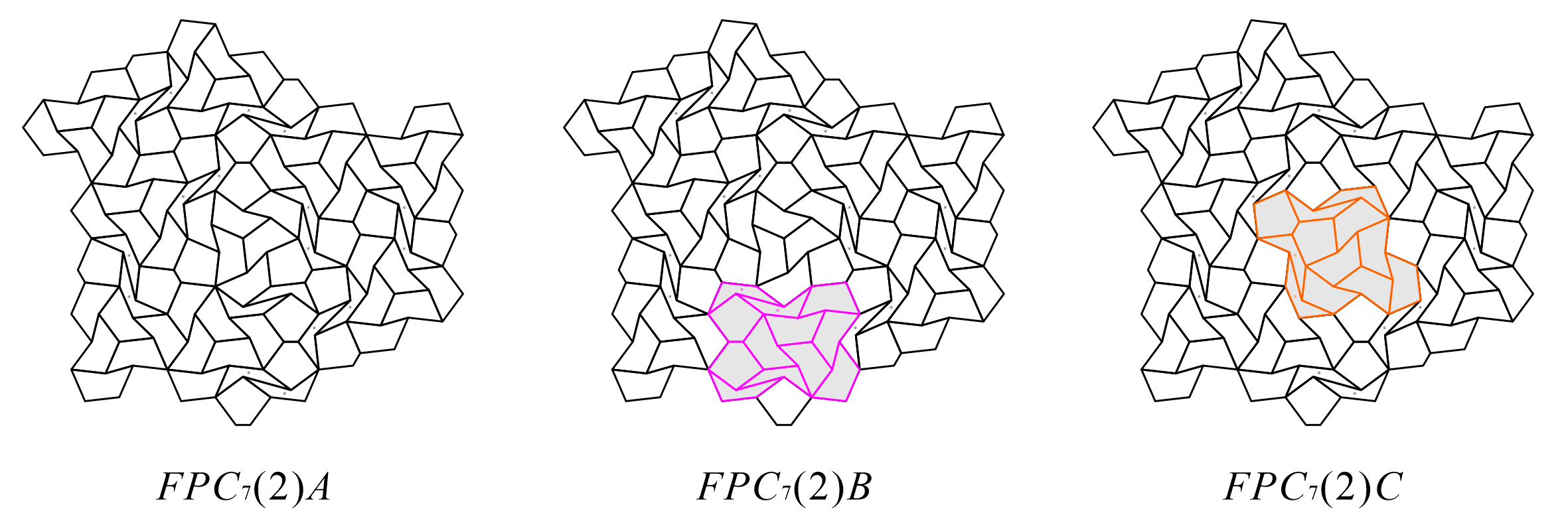} 
  \caption{{\small 
Locations where patterns of pentagons changed in 
$FPC_{7}(2)B$ and $FPC_{7}(2)C$ relative to $FPC_{7}(2)A$.
} 
\label{Fig.4-15}
}
\end{figure}

\begin{figure}[htbp]
 \centering\includegraphics[width=15cm,clip]{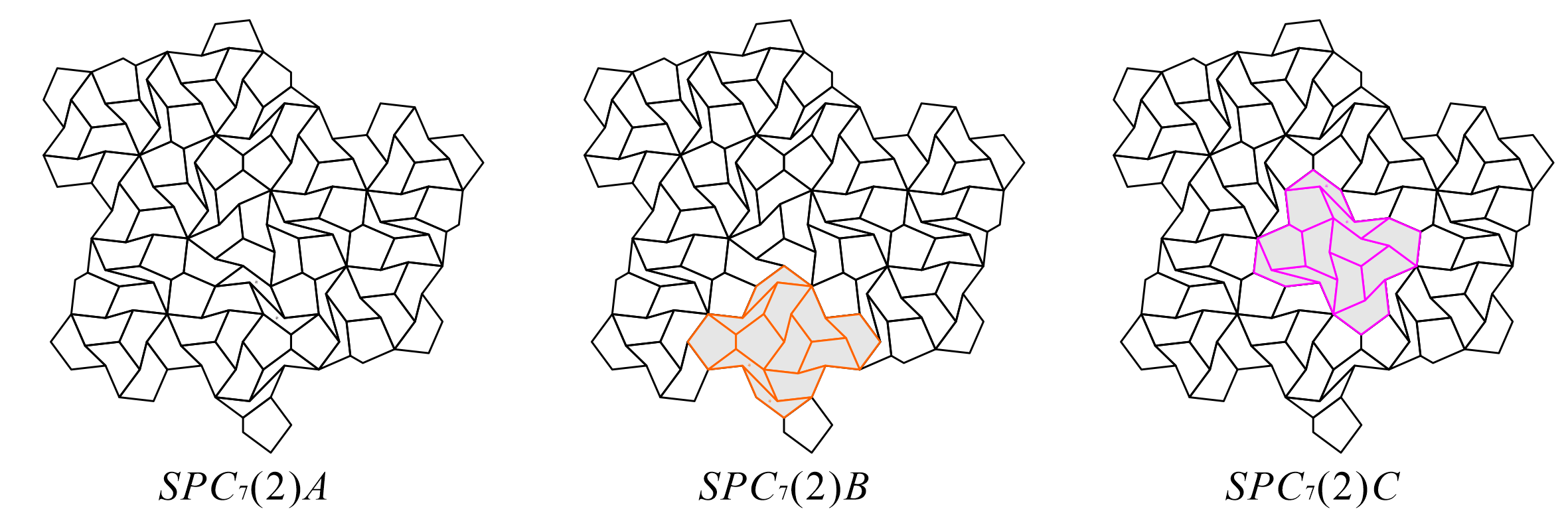} 
  \caption{{\small 
Locations where patterns of pentagons changed in 
$SPC_{7}(2)B$ and $SPC_{7}(2)C$ relative to $SPC_{7}(2)A$.
} 
\label{Fig.4-16}
}
\end{figure}

\begin{figure}[htbp]
 \centering\includegraphics[width=12cm,clip]{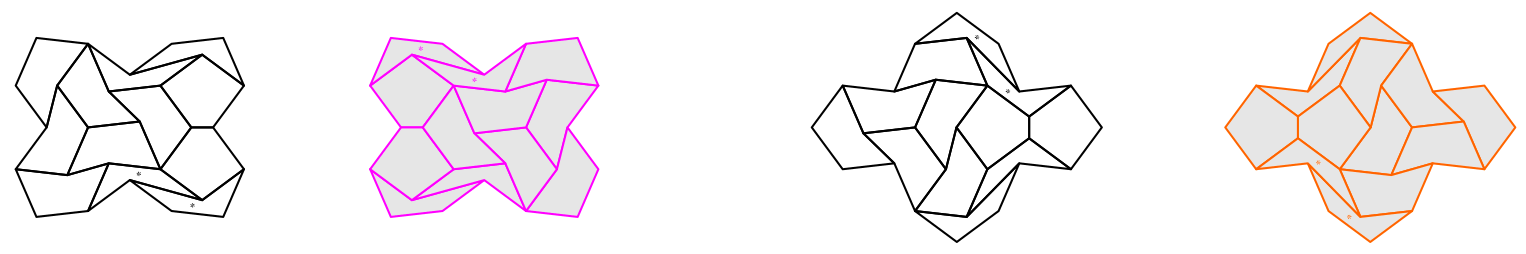} 
  \caption{{\small 
Components of change locations.
} 
\label{Fig.4-17}
}
\end{figure}

Figure~\ref{Fig.4-15} shows the locations where the pattern of pentagons changed in 
$FPC_{7}(2)B$ and $FPC_{7}(2)C$ relative to $FPC_{7}(2)A$. Similarly, Figure~\ref{Fig.4-16} 
shows the locations where the pattern of pentagons changed in $SPC_{7}(2)B$ and 
$SPC_{7}(2)C$ relative to $SPC_{7}(2)A$. As shown in Figures~\ref{Fig.4-15} and \ref{Fig.4-16}, 
two types of components exist at the locations where changes occur, as shown in 
Figure~\ref{Fig.4-17}. Both components exhibit $180^ \circ$ rotational symmetry in their outlines. 
However, they do not exhibit such symmetric property when considering the arrangement of 
the inner pentagons. As shown in Figure~\ref{Fig.4-10}, we can assume that Types B and C 
are created by replacing the Mystic regions in Type A with the Mystic formed by the reflected 
\mbox{Tile$(1, 1)$}. Meanwhile, we can assume that the clusters corresponding to Types B and C 
shown in Figures~\ref{Fig.4-15} and \ref{Fig.4-16} ($FPC_{7}(2)B$, $FPC_{7}(2)C$, $SPC_{7}(2)B$, 
and $SPC_{7}(2)C$) are created by replacing the clusters of Type A ($FPC_{7}(2)A$ and 
$SPC_{7}(2)A$) with the components corresponding to Figure~\ref{Fig.4-17}, rotated by $180^ \circ$.

\begin{figure}[tbp]
 \centering\includegraphics[width=15cm,clip]{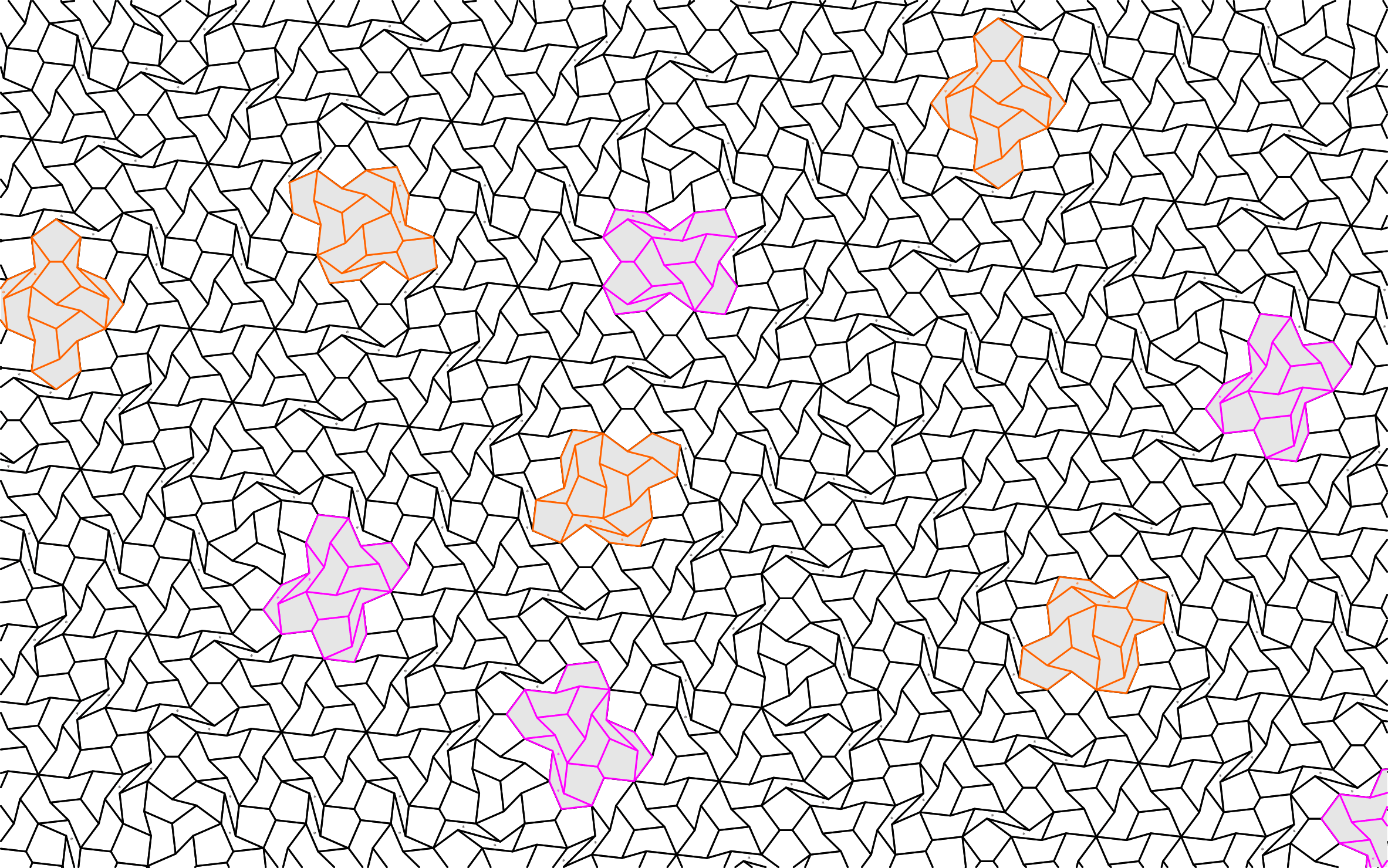} 
  \caption{{\small 
Tiling in which certain regions have been partially replaced by the two types of 
components shown in Figure~\ref{Fig.4-17}, based on the first pattern of $T_{s}$ 
shown in Figure~\ref{Fig.4-6}.
} 
\label{Fig.4-18}
}
\end{figure}

Figure~\ref{Fig.4-18} shows a tiling in which certain regions have been partially 
replaced by the two types of components shown in Figure~\ref{Fig.4-17}, based on the 
first pattern of $T_{s}$ shown in Figure~\ref{Fig.4-6} .

\renewcommand{\figurename}{{\small Figure}}

\renewcommand{\figurename}{{\small Figure}}

\renewcommand{\figurename}{{\small Figure}}


\renewcommand{\figurename}{{\small Figure}}


\section{Conclusions}
\label{section5}

In this study, the non-periodic tilings $T_{h}$ and $T_{s}$ with \mbox{Tile$(1, 1)$} were 
converted into pentagonal tilings without applying process (ii) in Section~\ref{section2}, 
``a tiling with rhombuses is converted into a tiling with pentagons after the original rhombus 
is divided such that it contains its own similar rhombuses (i.e., the original rhombus 
is subdivided into smaller similar rhombuses)." Based on the results of this study, 
the following conclusions were inferred:

\begin{itemize}
	\item \mbox{Tile$(1, 1)$} can generate countless patterns (design patterns created by 
	the arrangement of polygonal tiles) of tilings if the use of reflected tiles is allowed during 
	the tiling generation process.
	\item The non-periodic tilings generated using \mbox{Tile$(1, 1)$} based on $T_{h}$ 
	and $T_{s}$ (including the countless patterns of non-periodic tilings based on 
	$T_{s}$ shown in Section~\ref{subsection4.3}) can be converted into non-periodic 
	tilings with three types of rhombuses, as well as into non-periodic tilings with three 
	types of pentagons\footnote{ The periodic tiling shown in Figure~\ref{Fig.1-1}(b) 
	can be converted similarly (see \ref{appB}).}.
	\item When arbitrary replacements---namely, the arbitrary replacements of 
	locations corresponding to the regular hexagon and to the structure of Mystic in 
	$T_{s}$---are excluded, the non-periodic tilings with three types of pentagons 
	corresponding to $T_{h}$ and $T_{s}$ with \mbox{Tile$(1, 1)$} exhibit two distinct 
	patterns each. Furthermore, when arbitrary replacements of locations corresponding 
	to the structure of Mystic in $T_{s}$ are allowed, countless patterns of non-periodic 
	tilings with three types of pentagons corresponding to $T_{s}$ fall into two distinct series.
\end{itemize}

In connection with the above results, we have the following consideration:

\begin{itemize}
	\item In the non-periodic tiling $T_{s}$ consisting of rhombuses with acute angles of 
	$90^ \circ$, $60^ \circ$, and $30^ \circ$, we assume that special rhombuses exist 
	among the rhombuses with an acute angle of $30^ \circ$. Specifically, the rhombus 
	with an acute angle of $30^ \circ$ within \mbox{Tile$(1, 1)$}, which appears in the 
	overlap of two Mystic regions in $C_{7}(2)$ and $C_{6}(2)$, is assumed to be distinct 
	from the others.
\end{itemize}

In the (countless patterns of) pentagonal non-periodic tilings based on $T_{s}$ shown in 
Section~\ref{section4}, the pentagons corresponding to the rhombus with an acute angle 
of $30^ \circ$ must be used on both the anterior and posterior sides. Therefore, if the 
pentagons on the anterior and posterior sides are not congruent when superimposed, 
they can be regarded as different tiles. Thus, in the case where the original rhombus is not 
subdivided into smaller similar rhombuses, we can consider that the pentagonal non-periodic 
tilings based on $T_{s}$ require at least four types of tiles\footnote{ 
In the case where the original rhombus is not subdivided into smaller similar rhombuses, 
regardless of whether one distinguishes between the anterior and posterior side tiles, the 
pentagonal non-periodic tilings based on $T_{h}$ require at least three types of tiles.
}. According to Section 5.3.4 in \cite{Sugimoto_2022}, when $\theta = 60^ \circ$, a rhombus 
with an acute angle of $30^ \circ$ corresponds to a trapezoid, which is a degenerate 
pentagon. Because the trapezoid exhibits line symmetry, no distinction is made between 
its anterior and posterior sides. Therefore, in the case where the original rhombus is not 
subdivided into smaller similar rhombuses, when converted into pentagonal tiling with 
$\theta = 60^ \circ$ as described in Section 5.3.4 in \cite{Sugimoto_2022}, 
we consider that at least three types of tiles (a trapezoid and two types of 
convex pentagons) can form non-periodic tilings based on $T_{s}$, even if the 
tiles on the anterior and posterior sides are distinguished.

The illustrations of tilings for the pentagonal cases presented in Sections 5.3.2 to 5.3.7 
of \cite{Sugimoto_2022}, which are omitted in this manuscript due to space or capacity 
limitations, are available on the author's website \cite{Sugimoto_site_Ts} and 
\cite{Sugimoto_site_Th}.

Finally, the following questions remain open for future investigation:

\begin{itemize}
	\item Is it possible to assign a matching condition\footnote{ 
	Matching conditions specify how tiles must connect to form a valid tiling, which can sometimes 
	be represented by assigning colors or orientations to specific edges of the tiles.
	} to each rhombus to generate non-periodic 
	tilings with the three types of rhombuses corresponding to $T_{h}$ and $T_{s}$? 
	(If a matching condition can be established for $T_{s}$, we believe that rhombuses 
	with an acute angle of $30^ \circ$ are important).
	\item Is it possible to assign a matching condition to each pentagon to generate non-periodic 
	tilings with the three types of pentagons (Figures~\ref{Fig.3-4}, \ref{Fig.3-7}, \ref{Fig.4-6}, 
	and \ref{Fig.4-9}) corresponding to $T_{h}$ and $T_{s}$?
\end{itemize}



\bigskip
\appendix
\def\thesection{Appendix \Alph{section}}

\section{}
\label{appA}

\renewcommand{\thefigure}{\Alph{section}{-}\arabic{figure}}

\begin{figure}[htbp]
 \centering\includegraphics[width=15cm,clip]{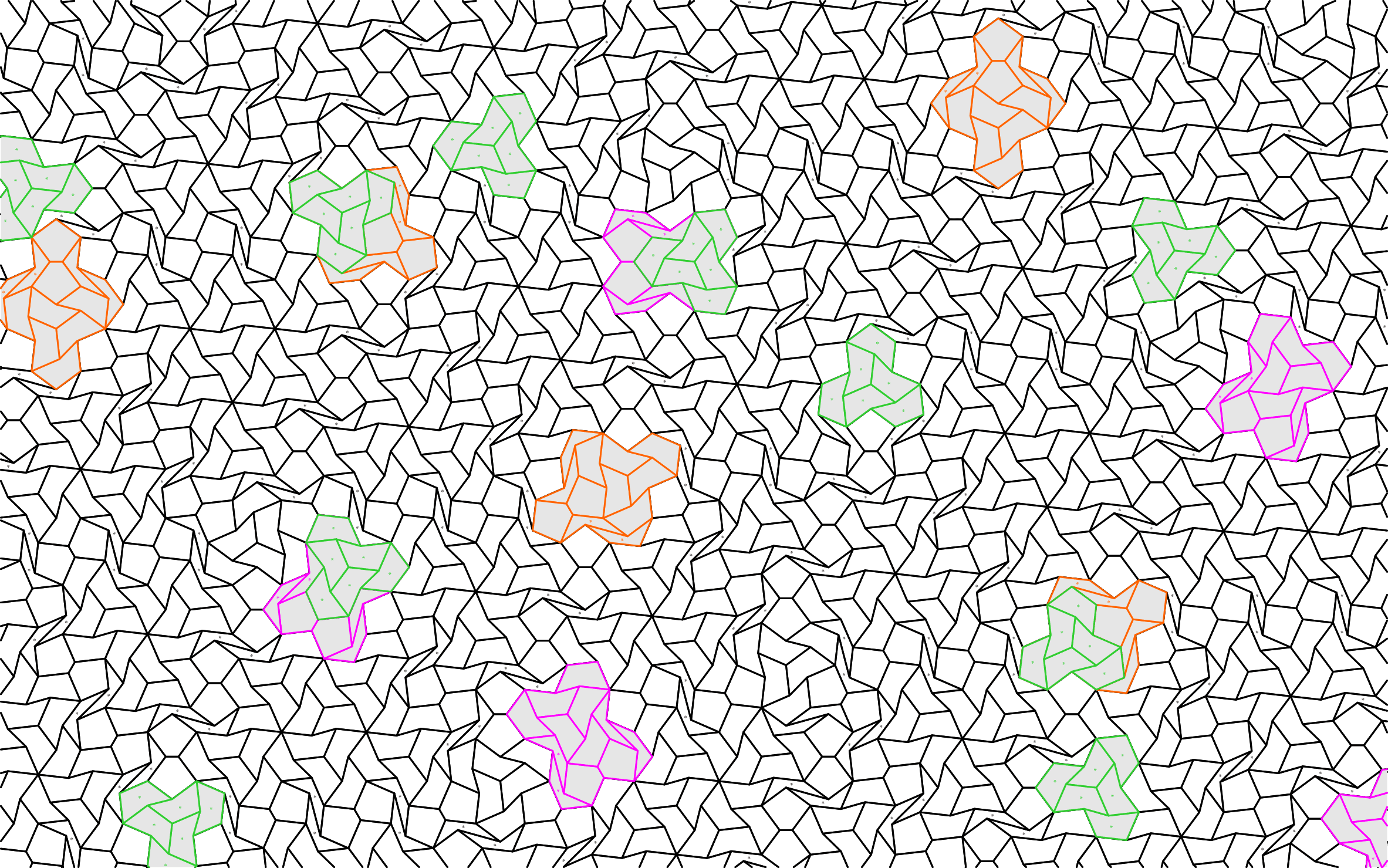} 
  \caption{{\small 
Tiling which was partially replaced arbitrarily, based on tiling in Figure~\ref{Fig.4-18}.
} 
\label{Fig.A-1}
}
\end{figure}

In this study, we did not arbitrarily replace locations corresponding to 
regular hexagons in the tilings during the conversion process. However, in 
this section, we briefly introduce a case in which locations corresponding 
to regular hexagons are arbitrarily replaced. Figure~\ref{Fig.A-1} shows a tiling in 
which a few locations corresponding to regular hexagons are replaced by 
units formed by three pairs of pentagons on the posterior side (marked 
with asterisks), based on the tiling shown in Figure~\ref{Fig.4-18}.

\renewcommand{\figurename}{{\small Figure}}

\begin{figure}[htbp]
 \centering\includegraphics[width=13cm,clip]{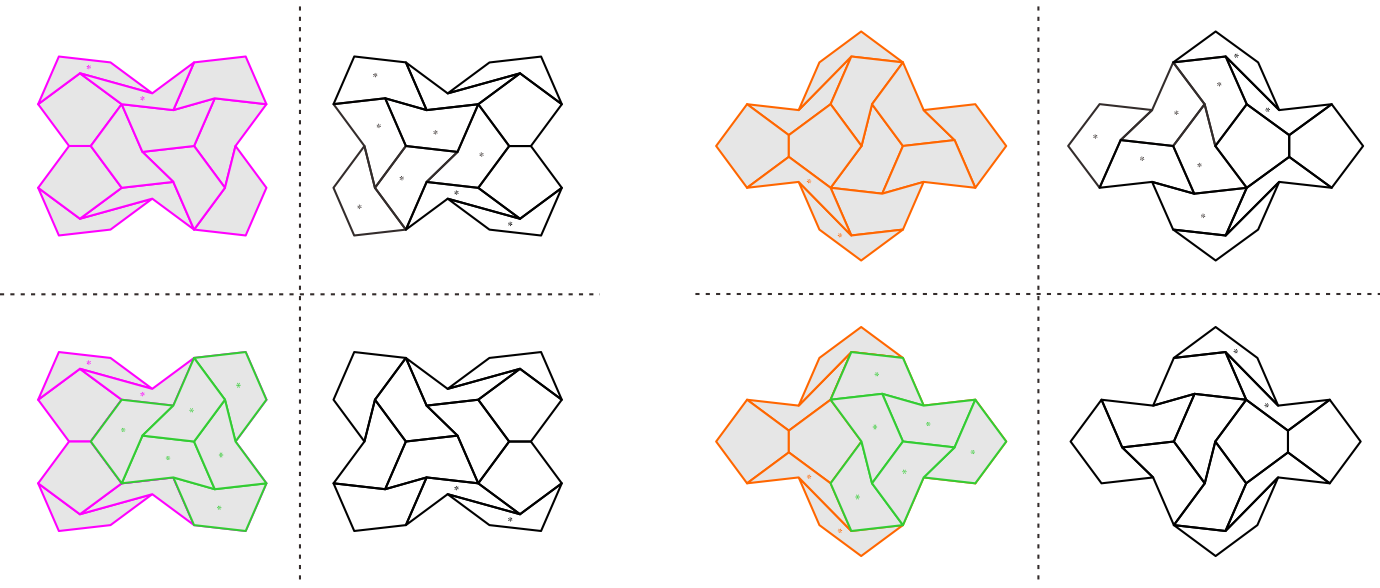} 
  \caption{{\small 
Properties of components in Figure~\ref{Fig.4-17}.
} 
\label{Fig.A-2}
}
\end{figure}

Because the outlines of the two types of components shown in Figure~\ref{Fig.4-17} 
exhibit line symmetry, their reflected components can be incorporated into 
the tiling. However, the reflected components shown in Figure~\ref{Fig.4-17} are 
identical to those obtained by replacing only the locations corresponding to 
regular hexagons in these components with a unit formed by three pairs of 
pentagons on the posterior side (see Figure~\ref{Fig.A-2}).

\renewcommand{\figurename}{{\small Figure}}


\section{}
\label{appB}

\begin{figure}[htbp]
 \centering\includegraphics[width=13.5cm,clip]{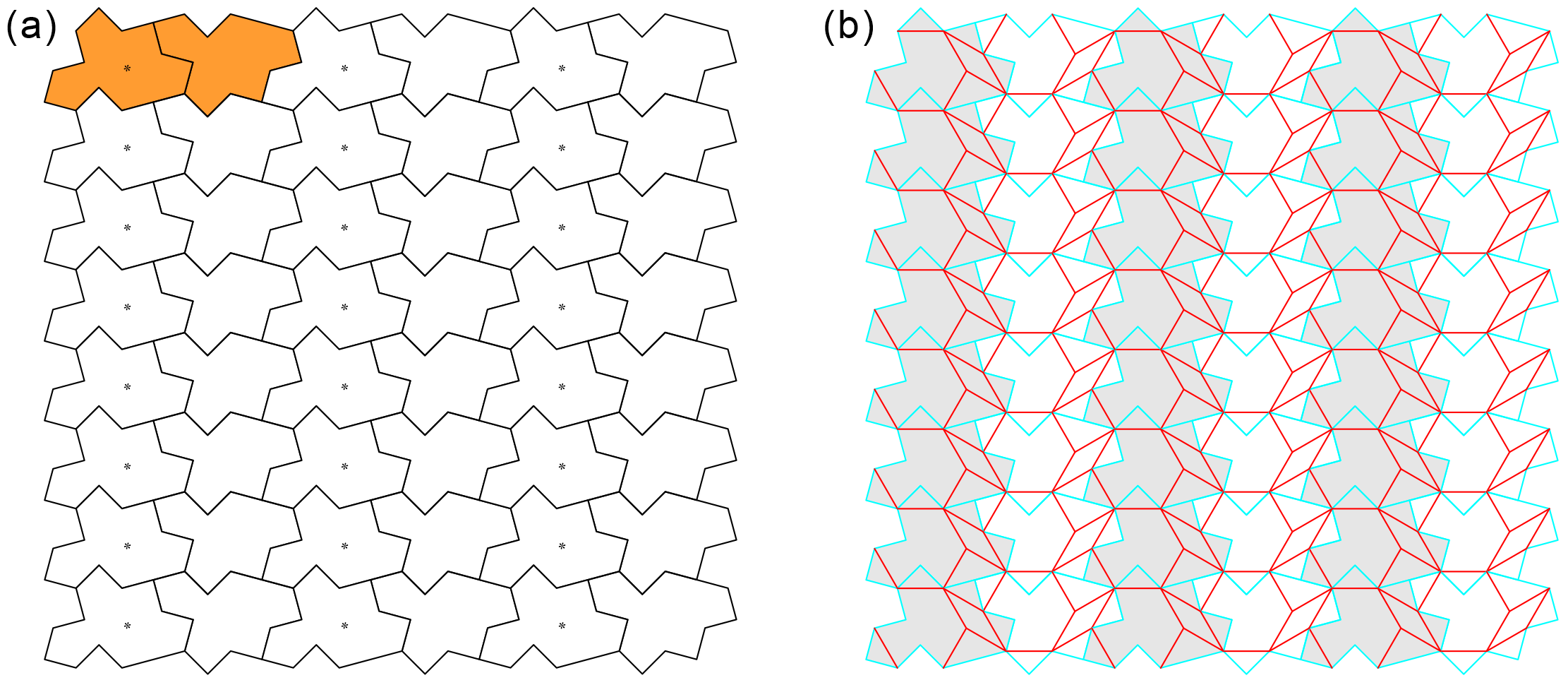} 
  \caption{{\small 
 Converting periodic tiling formed by \mbox{Tile$(1, 1)$} into tiling consisting of squares, 
 regular hexagons, and rhombuses with an acute angle of $30^ \circ$. The orange 
 region in the tiling is a translation unit.
} 
\label{Fig.B-1}
}
\end{figure}

\begin{figure}[htbp]
 \centering\includegraphics[width=13.5cm,clip]{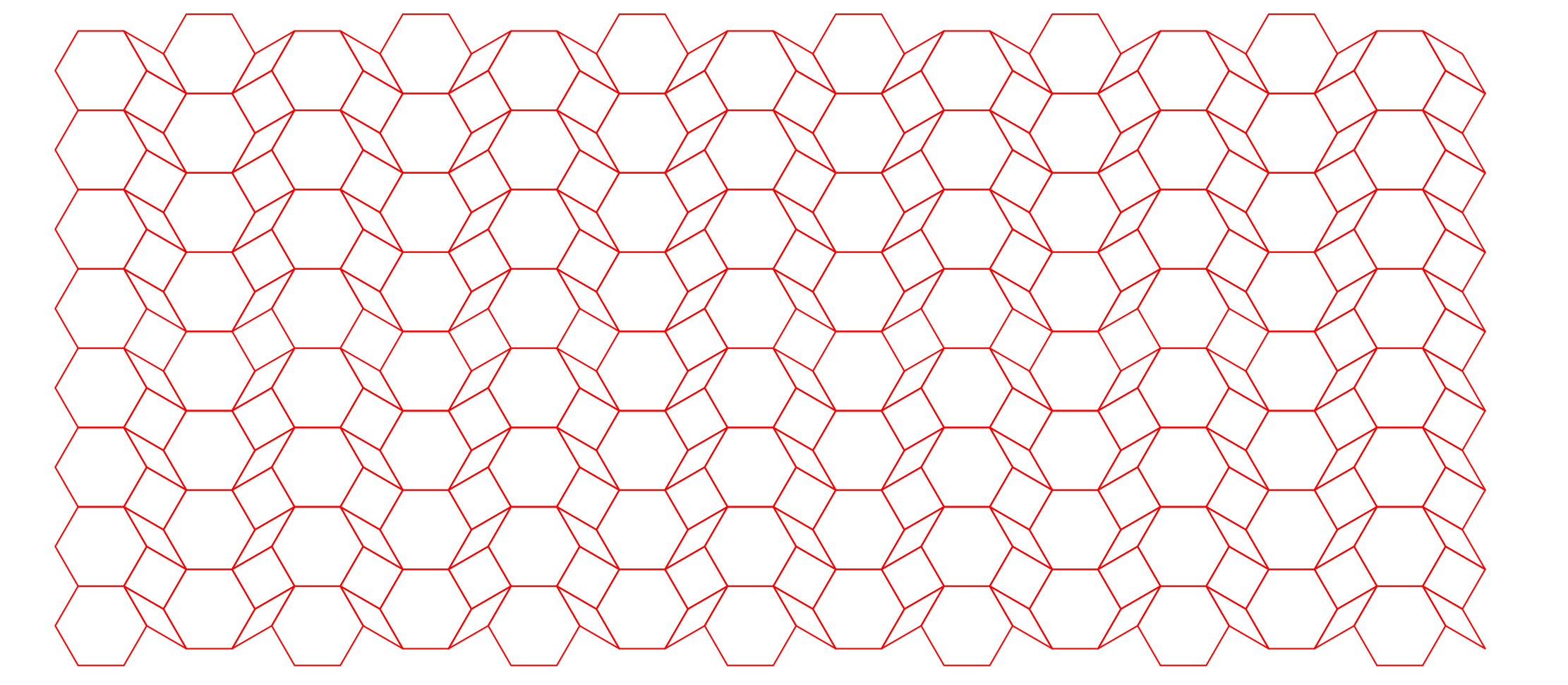} 
  \caption{{\small 
Periodic tiling formed by squares, regular hexagons, and rhombuses with an 
acute angle of $30^ \circ$.
} 
\label{Fig.B-2}
}
\end{figure}

If \mbox{Tile$(1, 1)$} is assigned a pattern incorporating decomposition lines, as shown 
in Figure~\ref{Fig.1-1}(c), then the periodic tiling formed by \mbox{Tile$(1, 1)$}, as shown 
in Figure~\ref{Fig.B-1}(a), can also be converted into a tiling consisting of squares, regular 
hexagons, and rhombuses with an acute angle of $30^ \circ$, as shown in Figure~\ref{Fig.B-1}(b). 
Now, let $T_{p }$ be the periodic tiling consisting of squares, regular hexagons, and rhombuses 
with acute angles of $30^ \circ$, created through this conversion, that is, a tiling exhibiting 
the pattern structure shown in Figure~\ref{Fig.B-2}.

\renewcommand{\figurename}{{\small Figure}}

\renewcommand{\figurename}{{\small Figure}}

\begin{figure}[t]
 \centering\includegraphics[width=13.5cm,clip]{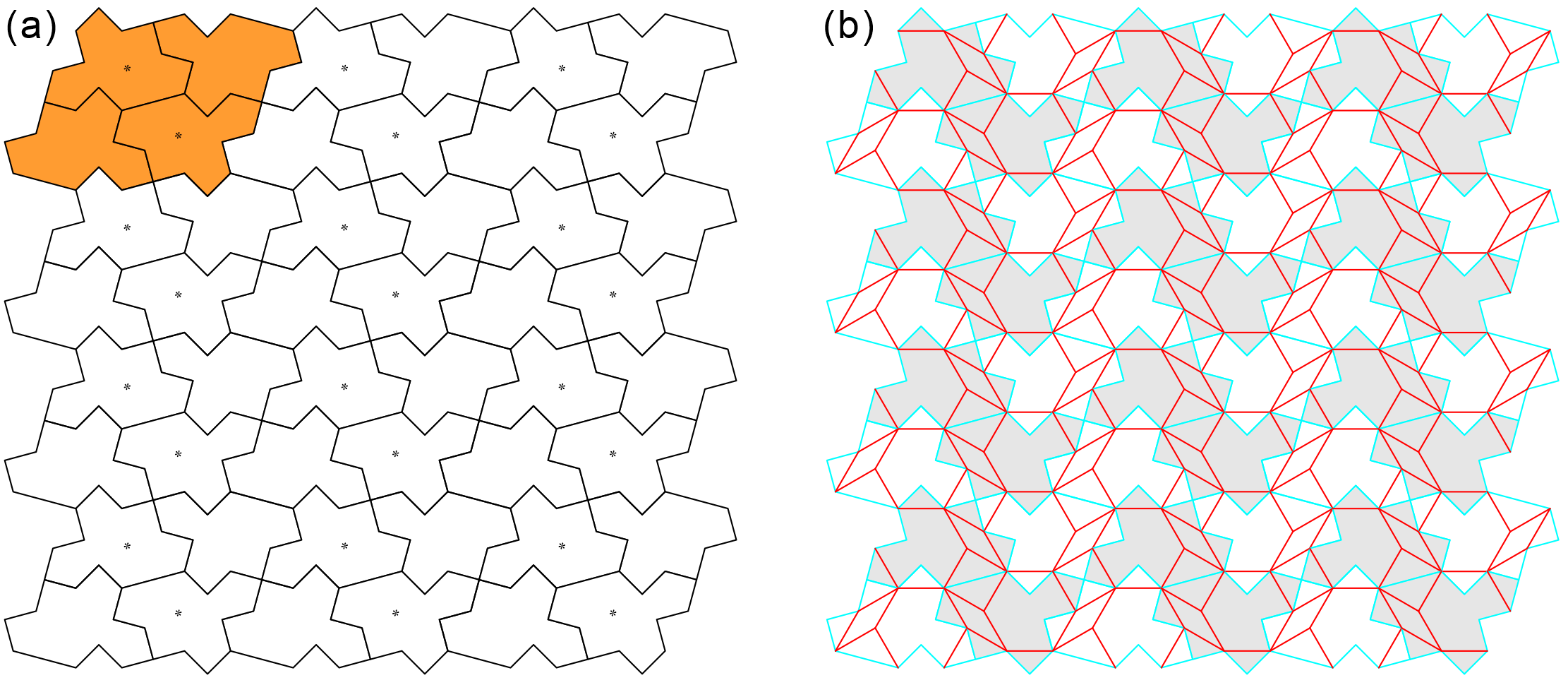} 
  \caption{{\small 
Periodic tiling generated by \mbox{Tile$(1, 1)$} other than tiling of 
Figure~\ref{Fig.B-1}. The orange region in the tiling is a translation unit.
} 
\label{Fig.B-3}
}
\end{figure}

\begin{figure}[t]
 \centering\includegraphics[width=13.5cm,clip]{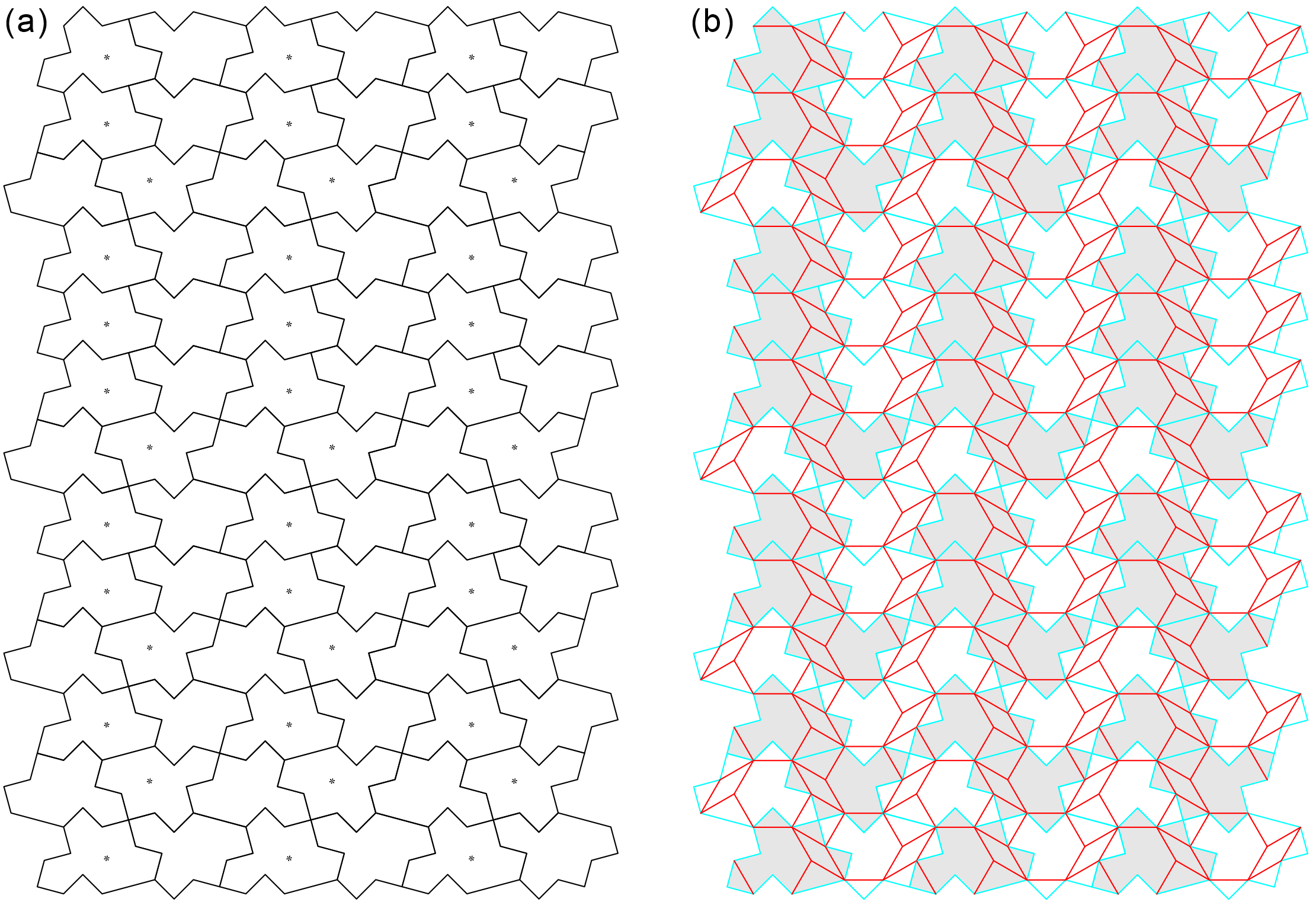} 
  \caption{{\small 
Example of tiling that can be converted into pattern in $T_{p }$ generated by \mbox{Tile$(1, 1)$}.
} 
\label{Fig.B-4}
}
\end{figure}

In addition to the tiling shown in Figure~\ref{Fig.B-1}, \mbox{Tile$(1, 1)$} can also generate 
a periodic tiling, as shown in Figure~\ref{Fig.B-3}(a). Furthermore, by focusing on the tilings 
in Figures~\ref{Fig.B-1}(a) and \ref{Fig.B-3}(a) and their translation units (units that can 
generate periodic tilings through translation alone), we observe that \mbox{Tile$(1, 1)$} can 
generate tilings, such as the one shown in Figure~\ref{Fig.B-4}(a), using belts consisting of these 
translation units. That is, the belts they form through horizontal translation in the same 
direction can be freely connected vertically. By applying the decomposition line pattern 
shown in Figure~\ref{Fig.1-1}(c) to \mbox{Tile$(1, 1)$} in the tilings in Figures~\ref{Fig.B-3}(a) 
and \ref{Fig.B-4}(a), we can convert them into periodic tilings consisting of squares, regular 
hexagons, and rhombuses with an acute angle of $30^ \circ$, as shown in Figures~\ref{Fig.B-3}(b) 
and \ref{Fig.B-4}(b), which are identified as $T_{p}$.

\renewcommand{\figurename}{{\small Figure}}

Therefore, the periodic tilings with \mbox{Tile$(1, 1)$} in Figures~\ref{Fig.B-1} and \ref{Fig.B-3} 
are both classified as $T_{p}$. Furthermore, a tiling in which belts are freely combined in the 
vertical direction by extending them parallel to the horizontal axis, as shown in Figure~\ref{Fig.B-4}, 
is also classified as $T_{p}$. Thus, \mbox{Tile$(1, 1)$} can generate a variety of patterns of 
tilings by combining bands of two types of translation units (i.e., there are a variety of patterns 
for $T_{p}$ with \mbox{Tile$(1, 1)$}). Because a regular hexagon can be divided into three 
rhombuses with an acute angle of $60^ \circ$, $T_{p}$ with \mbox{Tile$(1, 1)$} can be converted into 
non-periodic tilings consisting of rhombuses with acute angles of $90^ \circ$, $60^ \circ$, and 
$30^ \circ$. Thus, $T_{p}$ with \mbox{Tile$(1, 1)$} can be converted into periodic tilings 
with three types of pentagons.

\renewcommand{\figurename}{{\small Figure}}

\begin{figure}[htbp]
 \centering\includegraphics[width=14cm,clip]{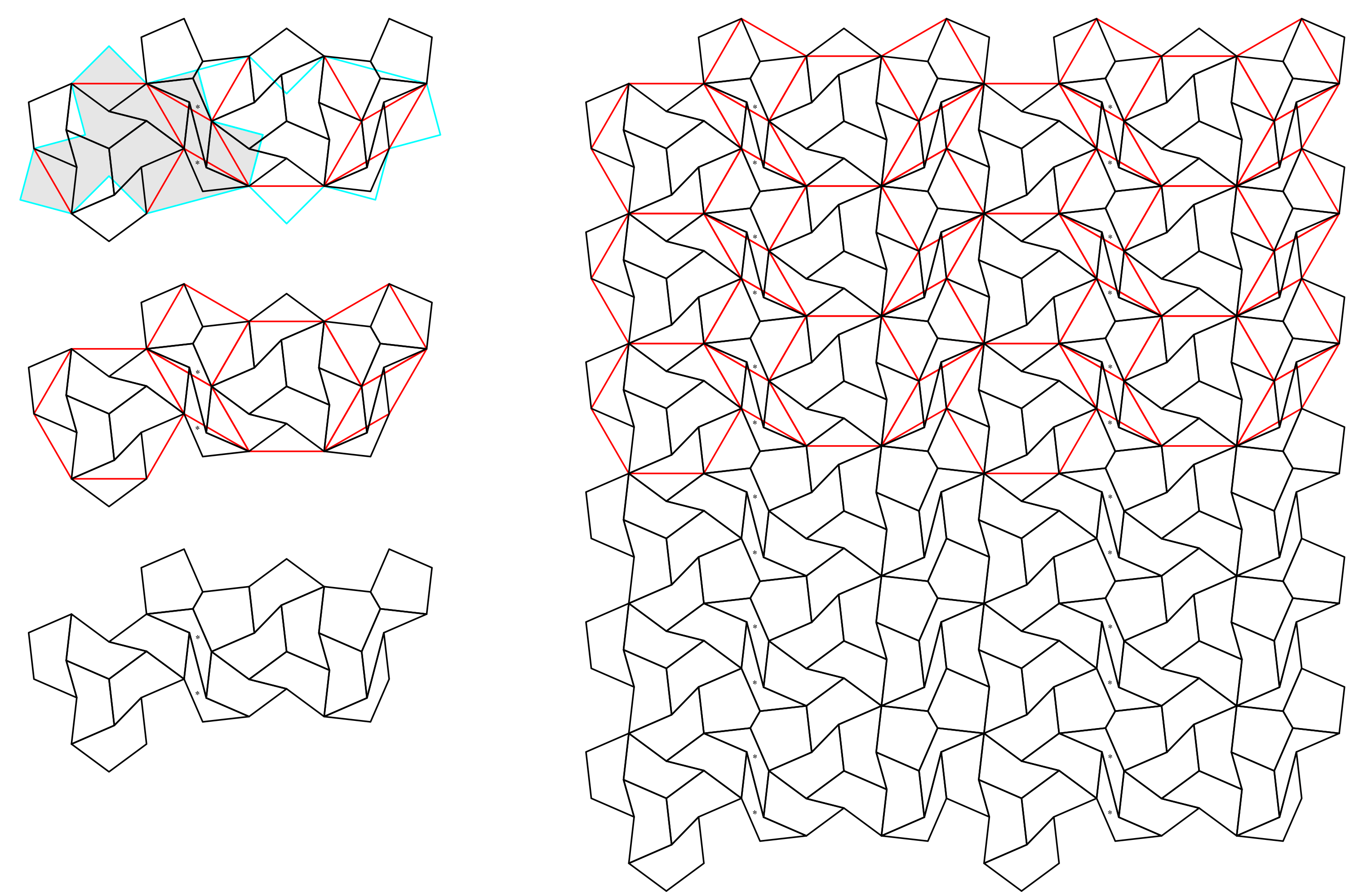} 
  \caption{{\small 
First pattern of $T_{p}$ with three types of pentagons (the first pattern when $T_{p}$ 
with \mbox{Tile$(1, 1)$} is converted into a tiling with three types of pentagons) in 
the case where the original rhombus is not subdivided into smaller similar rhombuses.
} 
\label{Fig.B-5}
}
\end{figure}

When $T_{p}$ with \mbox{Tile$(1, 1)$} is converted into tilings with three types of pentagons 
without subdivision of the original rhombus into smaller similar rhombuses (i.e., in the case where 
the original rhombus is not subdivided into smaller similar rhombuses), and the arbitrary 
replacements of locations corresponding to regular hexagons are excluded, the resulting periodic 
tilings with three types of pentagons exhibit the following two patterns:

Figure~\ref{Fig.B-5} shows the first pattern of $T_{p}$ (the first pattern in the case 
where $T_{p}$ with \mbox{Tile$(1, 1)$} is converted into a tiling with three types of pentagons) 
in the case where the original rhombus is not subdivided into smaller similar rhombuses. 
Specifically, it refers to the case in which the pentagons in the pair corresponding to the rhombus 
with an acute angle of $30^ \circ$, contained within \mbox{Tile$(1, 1)$} on the anterior side in 
Figure~\ref{Fig.1-1}, are marked without asterisks (the anterior side in Figure~\ref{Fig.2-2}), 
while those on the posterior side in Figure~\ref{Fig.1-1} are marked with 
asterisks (the posterior side in Figure~\ref{Fig.2-2}).

\begin{figure}[htbp]
 \centering\includegraphics[width=14cm,clip]{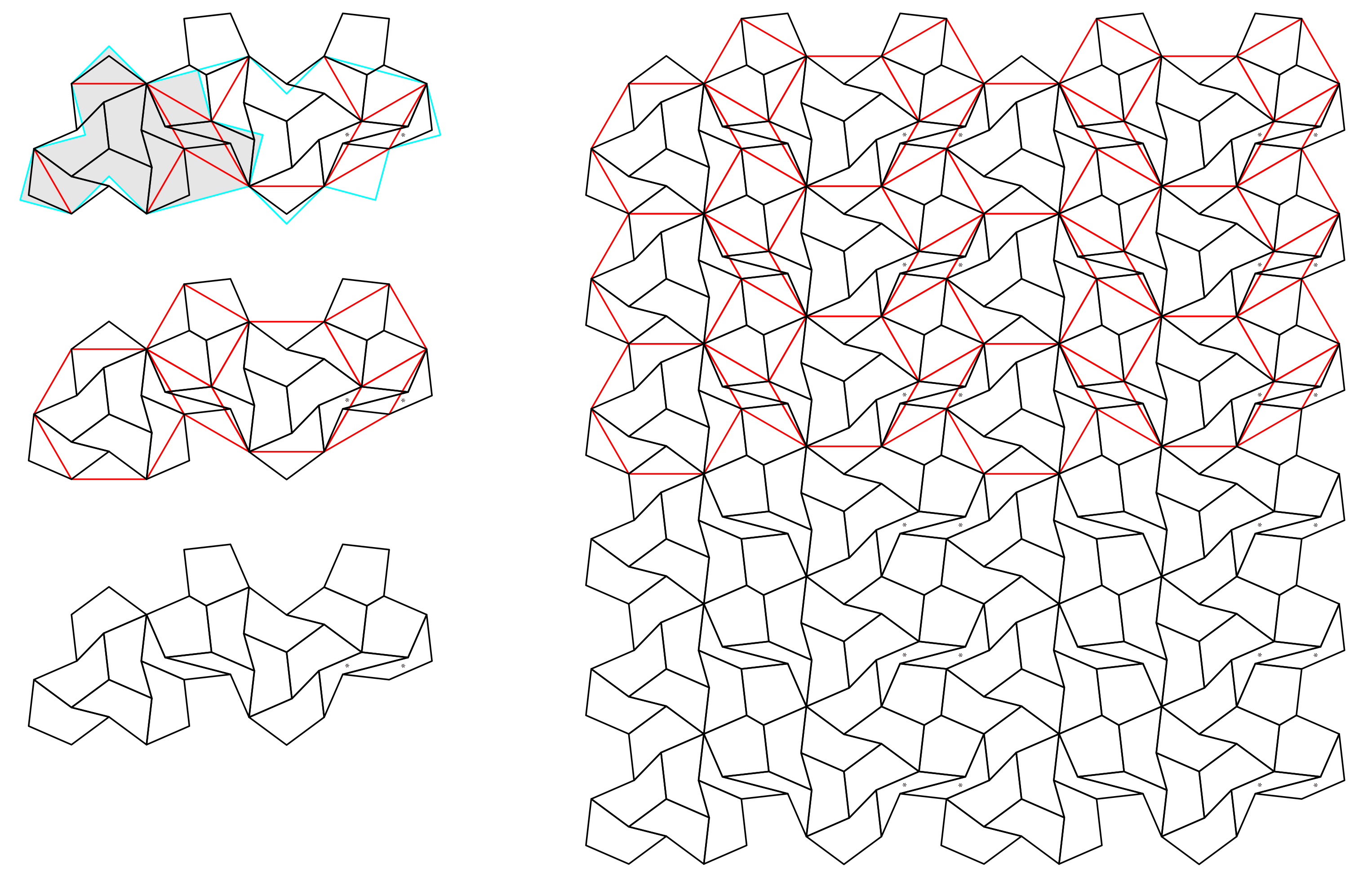} 
  \caption{{\small 
Second pattern of $T_{p}$ with three types of pentagons (the second pattern when $T_{p}$ 
with \mbox{Tile$(1, 1)$} is converted into a tiling with three types of pentagons) in the 
case where the original rhombus is not subdivided into smaller similar rhombuses.
} 
\label{Fig.B-6}
}
\end{figure}

Figure~\ref{Fig.B-6} shows the second pattern of $T_{p}$ (the second pattern in 
the case where $T_{p}$ with \mbox{Tile$(1, 1)$} is converted into a tiling with three types 
of pentagons) in the case where the original rhombus is not subdivided into smaller similar 
rhombuses. Specifically, it refers to the case in which the pentagons in the pair corresponding to 
the rhombus with an acute angle of $30^ \circ$, contained within \mbox{Tile$(1, 1)$} on the 
anterior side in Figure~\ref{Fig.1-1}, are marked with asterisks (the posterior side in 
Figure~\ref{Fig.2-2}), while those on the posterior side in Figure~\ref{Fig.1-1} are marked 
without asterisks (the anterior side in Figure~\ref{Fig.2-2}).

\renewcommand{\figurename}{{\small Figure}}

\renewcommand{\figurename}{{\small Figure}}

\end{document}